              \def\version{5 October 2020}	        	%
\documentclass[reqno,11pt]{amsart} 
\usepackage[T1]{fontenc}
\usepackage[utf8]{inputenc}
\usepackage{amsmath} 
\usepackage[mathscr]{eucal}
\usepackage{amssymb}
\usepackage{srcltx} 
\usepackage{dsfont}
\usepackage{hyperref}
\usepackage{color}
\usepackage{enumerate}
\usepackage{tikz}
\usetikzlibrary{arrows}
\usepackage{comment}
\usepackage{lscape}
\usepackage{graphicx}
\usepackage{bbm}

\numberwithin{equation}{section}
 
\def\one{\mathbbmss 1}
\def\emptyset{\varnothing} 
\def\lc{\lambda_{\rm c}}

\def\a{\alpha} 
\def\X{{\mathbf{X}}}
\def\Gg{{G_{\gamma}(\X^\la)}}
 
\def\g{\gamma} 
 
\def\d{{\rm d}} 
\def\dos{d_o}
\def\e{\varepsilon}

\def\lacc{\lambda^*_{\rm c}} 

\def\rG{r_{\rm B}} 

\def\la{\lambda}

\def\z{\mu} 
 
\def\L{\Lambda}

\newfam\Bbbfam 
\font\tenBbb=msbm10 
\font\sevenBbb=msbm7 
\font\fiveBbb=msbm5 
\textfont\Bbbfam=\tenBbb 
\scriptfont\Bbbfam=\sevenBbb 
\scriptscriptfont\Bbbfam=\fiveBbb 
\def\Bbb{\fam\Bbbfam \tenBbb} 

\def\NN{\mathrm{V}}
\def\NNN{\mathbf{V}}
\def\Psup{P_{\sup}}

\def\bs{\boldsymbol}
\def\No{{N_o}}

\def\SINR{\mathrm{SINR	}}

\newcommand{\R}     {\mathbb{R}} 
\newcommand{\Z}     {\mathbb{Z}} 
\newcommand{\N}     {\mathbb{N}} 
\renewcommand{\P}   {\mathbb{P}} 
 
\newcommand{\E}     {\mathbb{E}} 
\newcommand{\Q}     {\mathbb{Q}}

\def\1{{\mathchoice {1\mskip-4mu\mathrm l}      % Blackboard bold 1 
{1\mskip-4mu\mathrm l} 
{1\mskip-4.5mu\mathrm l} {1\mskip-5mu\mathrm l}}} 
 
\def\comment#1{} 
\newtheoremstyle{thm}{2ex}{2ex}{\itshape\rmfamily}{} 
{\bfseries\rmfamily}{}{1.7ex}{} 
 
\newtheoremstyle{rem}{1.3ex}{1.3ex}{\rmfamily}{} 
{\itshape\rmfamily}{}{1.5ex}{}

\renewcommand{\theequation}{\thesection.\arabic{equation}} 
 
\newtheorem{theorem}{Theorem}[section] 
\newtheorem{lem}[theorem]{Lemma} 
\newtheorem{prop}[theorem] {Proposition} 
\newtheorem{cor}[theorem]  {Corollary}

\newtheorem{step}{STEP} 
 
\theoremstyle{definition}
\newtheorem{defn}[theorem] {Definition} 
\newtheorem{thm}[theorem] {Theorem}

\newtheorem{remark}[theorem]  {Remark}

\renewcommand{\d}{{\rm d}} 
 
\newcommand{\eps}{\varepsilon}

\newcommand{\essinf}{{\operatorname {essinf}}}
\newcommand{\supp}{{\operatorname {supp}}} 
 
\newcommand{\dist}{{\operatorname {dist}}}

\DeclareMathOperator*{\esssup}{ess\,sup}

\def\supp{\mathrm{supp}}

\newcommand{\Ocal}   {{\mathcal O }}

\newcommand\numberthis{\addtocounter{equation}{1}\tag{\theequation}}
\renewcommand{\e}   {{\operatorname e }}

\definecolor{Red}{rgb}{1,0,0}

 \title[SINR percolation with random powers]{SINR percolation for Cox point processes\\  with random powers}
 \author[Benedikt Jahnel and András Tóbiás]{}
 
\begin{document}
 \maketitle
 \centerline{{\sc Benedikt Jahnel\footnote{Weierstrass Institute for Applied Analysis and Stochastics, Mohrenstra\ss e 39,
10117 Berlin, Germany \texttt{Jahnel@wias-berlin.de}} and András Tóbiás\footnote{Berlin Mathematical School, TU Berlin,
Straße des 17.~Juni 136, 10623 Berlin, Germany \linebreak
\texttt{Tobias@math.tu-berlin.de}}}}
\renewcommand{\thefootnote}{}
%\vspace{0.5cm}
%\centerline{\textit{(affiliations...)}}

% \centerline{{\sc Benedikt Jahnel\footnote{Weierstrass Institute for Applied Analysis and Stochastics, Mohrenstraße 39
%10117 Berlin, Germany
%\texttt{Benedikt.Jahnel@wias-berlin.de}} and András Tóbiás\footnote{
%Berlin Mathematical School, TU Berlin,
%Straße des 17.~Juni 136, 10623 Berlin, Germany
%\texttt{Tobias@math.tu-berlin.de}}}

\bigskip

\centerline{\small(\version)} 
\vspace{.5cm} 
 
\begin{quote} 
{\small {\bf Abstract:}} Signal-to-interference plus noise ratio (SINR) percolation is an infinite-range dependent variant of continuum percolation modeling connections in a telecommunication network. Unlike in earlier works, in the present paper the transmitted signal powers of the devices of the network are assumed random, i.i.d.~and possibly unbounded. Additionally, we assume that the devices form a stationary Cox point process, i.e., a Poisson point process with stationary random intensity measure, in two or higher dimensions. We present the following main results. First, under suitable moment conditions on the signal powers and the intensity measure, there is percolation in the SINR graph given that the device density is high and interferences are sufficiently reduced, but not vanishing. Second, if the interference cancellation factor $\gamma$ and the SINR threshold $\tau$ satisfy $\gamma \geq 1/(2\tau)$, then there is no percolation for any intensity parameter. Third, in the case of a Poisson point process with constant powers, for any intensity parameter that is supercritical for the underlying Gilbert graph, the SINR graph also percolates with some small but positive interference cancellation factor.
\end{quote}

%\vfill



\setcounter{tocdepth}{3}

%===========================================
%===========================================
%===========================================

\section{Introduction}\label{sec-Intro}
The study of percolation properties of infinite random graphs traces back many decades and many of the classical results are already available in textbook form, see for example the monographs~\cite{MR96,G99,BR06}. The first results for percolation in the continuum $\R^d$ were presented in the landmark paper by Gilbert~\cite{G61}, where non-trivial regimes of existence and absence of percolation (i.e., existence of an unbounded connected component) were established for the \emph{Poisson--Gilbert graph}, i.e., where the set of nodes is given by a homogeneous {\em Poisson point process} (PPP) and edges are drawn between two nodes whenever their distance is below a certain fixed positive connectivity threshold. 
%, consisting of vertices given by a homogeneous PPP $X^\la$ in $\R^2$ with intensity $\la>0$ and edges connecting any pair of vertices with distance less than $r>0$
The context of telecommunication systems was already mentioned there. 

In order to make the model more flexible, instead of using a fixed connectivity threshold, the nodes in the PPP can also be equipped with independent random radii, drawn from a common distribution, and any two nodes are connected by an edge whenever their distance is below the sum of the associated radii. In view of our topic in this manuscript, we call the resulting model a \emph{Poisson--Gilbert graph with random radii}, whose percolation properties can equivalently be expressed in terms of the corresponding \emph{Poisson--Boolean model with random radii}, and note that a wide range of results for percolation for this model are available, see for example~\cite{MR96,ATT16,G08,G09,GT18} and references therein. 

\medskip
In view of applications in wireless telecommunication systems, the extension of Poisson--Gilbert graphs towards Gilbert graphs based on {\em Cox point processes} (CPP), i.e., PPPs with random intensity measure, allows to study long-range communication properties in device-to-device networks where devices are placed according to a PPP in {\em random environment} that is represented by the intensity measure of the CPP. 
%Standard examples of asymptotically-essentially-connected environments with applications in telecommunications are Poisson--Voronoi, or Poisson--Delaunay tessellations, see for example~\cite{HJC17,CGH18}. 
Recently continuum percolation and further properties of such \emph{Cox--Gilbert graphs} were studied under certain conditions of stabilization and connectedness, see~\cite{HJC17,CGH18} and below. However, the edge-drawing rule remained, as in the classical case of the Poisson--Gilbert graph, based on a fixed connectivity threshold. 

In the very recent manuscript~\cite{JTC20}, first continuum percolation results are presented for Cox--Gilbert graphs where, as in the Poisson--Gilbert graph with random radii, each node is equipped with a random radius, and edges are placed between any two nodes whenever their distance is below the sum of the radii. In this case, again under certain stabilization and connectedness assumptions, most of the percolation properties of the Poisson--Gilbert graph with random radii can be reproduced also for this {\em Cox--Gilbert graph with random radii}. We note that here, similarly to \cite{MR96}, percolation properties of the Gilbert graphs are again expressed in terms of the underlying \emph{Cox--Boolean models with random radii}.

\medskip
Moving beyond a setting where edges are placed between pairs of points based on their mutual distance and their associated radii, another line of research aims towards a different kind of extension of Gilbert graphs with respect to the edges. Starting with the papers~\cite{DBT05,DF06}, still based on a homogeneous PPP in $\R^2$, the edge-drawing mechanism is replaced by a highly non-local rule using the {\em signal-to-interference ratio} (SINR), which we describe precisely and more generally in~\eqref{SINR}. In words, very roughly, a pair of Poisson points is connected by an edge only if the weighted distance between the points is sufficiently small compared to the accumulated weighted distances of all the other points, the so-called {\em interference}. This definition is very much inspired by applications in wireless telecommunication networks, where the success of a transmission between network components is highly dependent on the relative signal strength between the components compared to the other (unwanted) signals present in the medium. In the simplest case, only the relative distances between points enter the SINR, giving rise to the SINR graph on PPPs, or the~\emph{Poisson-SINR graph}. Let us note that this definition introduces long-range, or even infinite-range dependencies for the construction of edges into the system. This is the setting considered in~\cite{DBT05,DF06}, where, using comparison techniques with the Poisson--Gilbert graph, again non-trivial percolation properties could be established. Let us mention that the SINR graph has very different monotonicity properties compared to the Poisson--Gilbert graph, which makes it more interesting but also harder to analyze. To see this, note that due to the definition of the SINR, an increase in the number of points also leads to an increase in the unwanted interference and thus to the potential loss of edges. On the other hand, for the Poisson--Gilbert graph, the connectivity increases if points are added into the system. 

Now, similar to the generalization of a Poisson--Gilbert graph towards a Poisson--Gilbert graph with random radii, also SINR graphs can be generalized in the sense that each point is equipped with an independent power random variable. These powers enter the definition of the SINR as presented in~\eqref{SINR} and represent the individual signal strengths of the network components. For the case of the SINR graph with random powers based on PPPs, or the \emph{Poisson--SINR graph with random powers}, the paper~\cite{KY07} presents first results similar to the assertions presented in~\cite{DBT05,DF06} under very strong boundedness assumptions on the powers. 
% under very strong boundedness assumptions on the fadings (in particular, boundedness from above and away from 0) and obtained similar assertions as \cite{DBT05,DF06} in the constant-power case. 
Let us note that the definition of an SINR graph with random powers already occurs in \cite{DBT05}, but the only proven result of this paper for this setting is about degree bounds (cf.~Section~\ref{sec-strategy2}).
The first steps towards understanding the case of unbounded powers in the Poisson--SINR graph with random powers were made recently in~\cite{L19}.
%in the recent master's thesis of R.~L\" offler at TU Berlin \cite{L19}, co-advised by the second and the third author of the present paper. 
In this master's thesis, supervised by the authors, results in the spirit of one of our main theorems, Theorem~\ref{thm-randomfadings}, are presented under much stronger assumptions and only for the case of an underlying homogeneous PPP.
%
%
%it was shown that for $d \geq 2$, finiteness of $\lambda_{\z}$ holds under than presented in Theorem~\ref{thm-randomfadings}. 
%on $\ell$ that are analogous to Assumption $(\ell)$ in case (a) $P_0$ is integrable and $\ell$ has bounded support or (b) $P_0$ has some exponential moments. \color{red} Regine claims that the critical intensity is the same as without interference, but I don't think that this follows from her proofs. I will discuss it with her and ask her to change it in case if it's really not correct. \color{black} The approach of the proof of this statement is similar to the main statement of \cite{DF06}, however, it additionally uses discrete percolation techniques originating from \cite{HJC17,T18}. 
\cite{L19} also provides sufficient conditions for the absence of percolation for small intensities of the PPP. Before \cite{L19}, no positive assertions about percolation in an SINR graph with unbounded random powers had been known in the literature; regarding the case of bounded random powers, see also Section~\ref{sec-StratProof}. 

\medskip
On the other hand, in~\cite{T18}, the two extensions described above were for the first time considered jointly, giving rise to the SINR graph based on CPPs, the \emph{Cox-SINR graph}, however without random powers. There it was established that for sufficiently large intensities and sufficiently connected environments, the Cox-SINR graph percolates almost surely at least for non-vanishing interference. In~\cite[Section 4.2.3.4]{T19b} it was anticipated that the case of random but bounded powers can easily be handled via the same methods, see Section~\ref{sec-StratProof1} for further details.

\medskip
The present manuscript now completes this line of research by analyzing the {\em Cox-SINR graph with random powers}, which are also not necessarily bounded. More precisely, in our first main result, Theorem~\ref{thm-randomfadings}, we present sufficient conditions for the existence of a supercritical percolation phase, i.e., a non-trivial regime for the intensity of the underlying CPP and non-vanishing interference, such that the Cox-SINR graph with random powers percolates. This substantially extends the results of \cite{DBT05, DF06, FM07} from the case of a homogeneous PPP in $\R^2$ with constant powers to the one of a CPP in $\R^d$, $d \geq 2$, with random and possibly unbounded powers, combining the methods of \cite{T18} for the case of a CPP with constant powers and the ones of \cite{L19} about the case of a PPP with random powers (both in dimensions 2 or higher). We will discuss the relation of Theorem~\ref{thm-randomfadings} to these results in detail in Section~\ref{sec-StratProof}. 
%The proof of Theorem~\ref{thm-randomfadings} uses some arguments of the proof of a previous result, Proposition~\ref{prop-T18}, which covers the case of bounded powers but does not tell anything about the case $P_{\sup}=\infty$.

In our second main result, Theorem~\ref{thm-2degree1cluster}, we establish a uniform bound on the strength of the interference, above which no percolation is possible. In essence, this theorem claims that there is no percolation in the SINR graph whenever degrees of its vertices are bounded by 2. The fact that SINR graphs with non-vanishing interference have bounded degrees originates from \cite[Theorem 1]{DBT05}; however, in that paper, only the trivial assertion that SINR graphs with degrees bounded by 1 do not percolate was proven, and this has not been improved before the present paper.

Finally, in our third main result, Theorem~\ref{thm-PPPlambdaequality}, we state that in the case of the Poisson-SINR graph with constant powers, indeed, the critical intensity for percolation in the presence of interference can be represented via the critical intensity of an associated Poisson--Gilbert graph. This result extends the two-dimensional statement \cite[Theorem 1]{DF06} to higher dimensions, whereas its proof is rather different from the one in \cite{DF06}.

\medskip
The organization of the manuscript is as follows. In Section~\ref{sec-IntroMainRes}, we present the setting and main results but postpone the introduction of our main technical conditions for the CPP, namely stabilization and asymptotic essential connectedness. These conditions are presented in Section~\ref{sec-Ex} together with examples for CPPs for which our main results are applicable. In Section~\ref{sec-StratProof} we present the proof strategies for our main results and give further background on how they relate to previous results in the literature. Finally, in Section~\ref{sec-Proofs}, we give the detailed proofs.

%===========================================
%===========================================
%===========================================

\section{Setting and main results}\label{sec-IntroMainRes}
\begin{figure}[!htpb]
\centering
\input{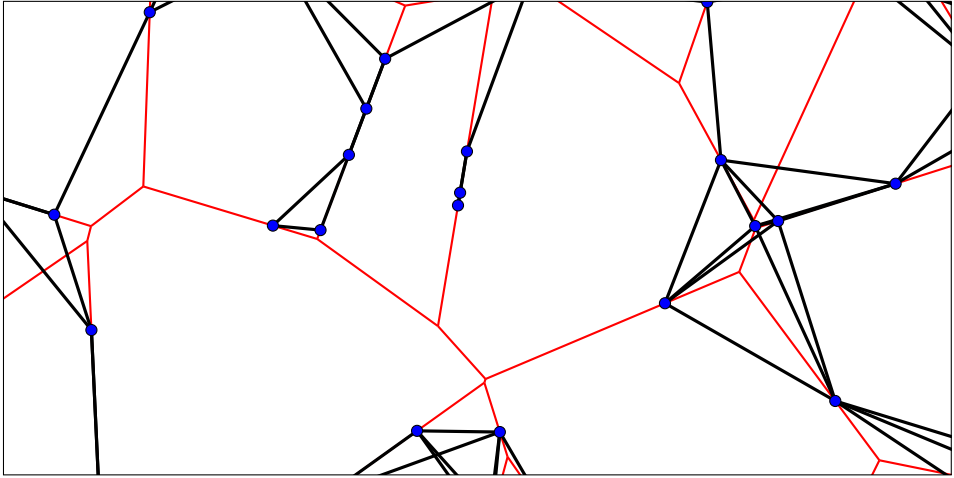}
\caption{A typical realization of a Cox-SINR graph (with blue vertices and black edges) with directing measure given by the edge-length measure of a two-dimensional Poisson--Voronoi tessellation (in red) in a box, with $\No=1$, $\g=0.3$, $\tau=0.2$, constant powers equal to $1$, and a suitable path-loss function $\ell$. 
}\label{fig-SINRGraph}
\end{figure}

For $\lambda>0$, let $\X^\la=\{(X_i,P_i)\}_{i \in I}$ be an i.i.d.~marked \emph{Cox point process} (CPP) in $\R^d\times[0,\infty)$ for $d \geq 1$, with directing measure $\la\L\otimes\z$ where $\L$ is stationary with $\E[\L(Q_1)]=1$ and $Q_n=[-n/2,n/2]^d$ for $n >0$. Here, $\mu$ is a Borel probability measure on $[0,\infty)$, the common distribution of the marks $P=\{P_i\}_{i\in I}$. We consider the \emph{SINR graph} with vertex set given by the first component of $\X^\la$, which we denote by $X^\la=\{X_i\}_{i\in I}$. Here, every pair $X_i, X_j\in X^\la$ of vertices with $X_i\neq X_j$, is connected by an edge if and only if
 \begin{equation}\label{SINR}
 \begin{split}
 P_i\ell(|X_i-X_j|)&>\tau\big(\No + \g \sum_{k \in I \setminus \{ i,j\}} P_k \ell(|X_k-X_j|)\big) \quad\text{and }\\
 P_j\ell(|X_i-X_j|)&>\tau\big(\No + \g \sum_{k \in I \setminus \{ i,j\}} P_k \ell(|X_k-X_i|)\big).
 \end{split}
 \end{equation}
In \eqref{SINR}, $\tau>0$ is fixed and called the \emph{SINR threshold}, the constant $\No\ge 0$ represents \emph{noise}, $r\mapsto\ell(r)\in [0,\infty)$ is referred to as the \emph{path-loss function} and $\g\ge0$ is called the \emph{interference-cancellation factor}. The random variables $P$ are often called random \emph{powers} and the term 
$$I(X_i,X_j, \X^\la)=\sum_{k \in I \setminus \{ i,j\}} P_k \ell(|X_k-X_j|)$$ 
is referred to as \emph{interference}.
We will use the notation $\Gg$ to indicate the SINR graph, suppressing the dependencies on $\tau$, $\No$ and $\ell$, but highlighting the dependence on $\g$, see also Figure~\ref{fig-SINRGraph} for an illustration.  We refer to \cite[Section 1]{DBT05} for further interpretation of the modeling parameters.

\medskip
The SINR graph has a nice interpretation in the study of device-to-device telecommunication systems where the devices $X^\la$ can communicate directly with each other if their mutual distance, represented by the path-loss function, and their individual powers, are sufficiently strong to overcome thermal noise plus all the interference coming from the other devices. If this is the case, then the possibility to communicate is represented by an undirected edge. 
%The SINR graph has been the subject of, by now, a large body of works, which we will further elaborate on in Section~\ref{sec-StratProof}. 

\medskip
Our main interest lies in percolation properties of the SINR graph, as has been first studied in~\cite{DBT05, DF06, FM07}. A \emph{cluster} in $\Gg$ is a maximal connected component. We say that $\Gg$ \emph{percolates} if $\Gg$ contains an unbounded component. Clusters and percolation are defined analogously in the case of any graph having a vertex set that is included in $\R^d$, $d \geq 1$, and is locally finite. Here we focus on the following key quantities. First, the \emph{critical interference-cancellation factor} is defined as 
 \[ \gamma(\lambda) = \sup \big\{ \gamma>0 \colon \P(\Gg\text{ percolates} ) >0 \big\}. \numberthis\label{gammastarlambda} \]
In words, it represents the maximal amount of interference that can be added to the system and still maintain percolation. Second, the \emph{critical intensity} is defined as 
 \[ \lacc=\inf \{ \lambda>0 \colon \gamma(\lambda')>0,~\forall \lambda'>\lambda \}, \numberthis\label{lambdacritSINR} \]
which describes the smallest intensity such that for all larger intensities the addition of a small amount of interference does not destroy percolation. 
%Let us mention that SINR graphs can be seen as subgraphs of Gilbert graphs

\medskip
For the statement of our first main result, we assume certain decorrelation and connectivity properties for the directing measure $\L$ of the underlying CPP. The precise definitions for $\L$ to be \emph{stabilizing}, \emph{$b$-dependent}  or \emph{asymptotically essentially connected} are technical and will be presented in Definitions~\ref{defn-stabilization} and~\ref{defn-asessconn} in Section~\ref{sec-Ex}, where we will also mention a number of relevant examples of random measures satisfying these definitions. We denote by $P_o$ a generic power random variable distributed according to $\z$ and put $\Psup=\esssup \mu$. We call the pathloss function $\ell$ {\em well-behaved} if $\ell$ is continuous, constant on $[0,\dos]$ for some $\dos \geq 0$, strictly decreasing on $[\dos,\infty) \cap \supp(\ell)$, and satisfies $\int_0^\infty r^{d-1} \ell(r) \d r <\infty$. Our first result establishes percolation for the SINR graph based on CPPs with random powers. 
\begin{thm}\label{thm-randomfadings}
Let $d \geq 2$, $\No,\tau>0$, $\Psup=\infty$, $\Lambda$ be stabilizing and $\ell$ be well-behaved. Then, $\lacc<\infty$ holds if at least one of the following conditions is satisfied:
\begin{enumerate}
\item\label{second-fading} 
$\ell$ has unbounded support, $\Lambda$ is $b$-dependent for some $b>0$, and $\E[\exp(\a \Lambda(Q_1))]<\infty$ as well as $\E[\exp(\a P_o)]<\infty$ holds for some $\a>0$, or
\item\label{first-fading} $\ell$ has bounded support, $\E[P_o]<\infty$, and 
$\Lambda$ is asymptotically essentially connected, or
\item\label{third-fading} $\ell$ has bounded support, $\E[P_o]<\infty$, and  $\sup \supp(\ell)$ is larger than $c$, where $c$ is a finite constant depending on $\Lambda,\tau$ and $\No$. 
\end{enumerate}
\end{thm}
As discussed in the introduction, Theorem~\ref{thm-randomfadings} extends similar results known for the case of a homogeneous PPP in $\R^2$ with constant powers to the one of a CPP in $\R^d$, $d \geq 2$ with random and possibly unbounded powers.
%, combining the methods of \cite{T18} for the case of a CPP with constant powers and the ones of \cite{L19} about the case of a PPP with random powers (both in dimensions 2 or higher). We will discuss the relation of Theorem~\ref{thm-randomfadings} to these results in detail in Section~\ref{sec-StratProof}. Before the master's thesis \cite{L19}, no positive assertions about percolation in an SINR graph with unbounded random powers had been known in the literature; regarding the case of bounded random powers, see also Section~\ref{sec-StratProof}.
%%The proof of Theorem~\ref{thm-randomfadings} uses some arguments of the proof of a previous result, Proposition~\ref{prop-T18}, which covers the case of bounded powers but does not tell anything about the case $P_{\sup}=\infty$. 
%Condition \eqref{third-fading} in Theorem~\ref{thm-randomfadings} can easily be omitted as long as we still assume that $\ell$ is continuous on $[0,\infty)$; nevertheless, we assume that it holds because it is the usual assumption on SINR percolation \cite{DBT05,DF06} and it makes the proof less technical. \color{red} Statement about $\lambda_{\z}>0$? \color{black} 
%Let us note that Condition~\eqref{third-fading} in Theorem~\ref{thm-randomfadings} can easily be omitted as long as we still assume that $\ell$ is continuous on $[0,\infty)$; nevertheless, we will use Condition~\eqref{third-fading} because it is the usual assumption on SINR percolation in \cite{DBT05,DF06} and it makes the proof less technical. 
Let us note that the complementary assertion that $\lacc>0$ can be deduced in certain cases based on recent results on Cox--Gilbert graphs with random radii, see~\cite{JTC20} and Section~\ref{sec-StratProof1} for more details.

\medskip
Our second main result establishes a uniform upper bound on the critical interferen\-ce-cancellation factor. For this assertion we assume the basic non-degeneracy property that $X^\lambda$ is \emph{nonequidistant}, which is satisfied for a very large class of CPPs including many relevant examples for wireless telecommunication systems, see Section~\ref{sec-Ex}. This means that for all $i,j,k,l\in I$, $| X_i-X_j|=| X_k-X_l|>0$ implies $\{ i,j\} = \{ k,l \}$ and $| X_i|=| X_j|$ implies $i\neq j$, almost surely.  Clearly, this property implies that the point process $X^\lambda$ is simple, further, if nonequidistantness of $X^\lambda$ holds for some $\lambda>0$, then it holds for all $\lambda>0$. As for a (pathological) counterexample, note that for $\Lambda$ being the sum of Dirac measures at the points of the randomly shifted lattice $\Z^d+U$, where $U$ is a uniform random variable in $[0,1]^d$, the associated CPP is simple, stationary, but not nonequidistant.
\begin{thm}\label{thm-2degree1cluster}
Let  $d \geq 1$, $\No\geq 0$ and $\tau,\lambda>0$, and assume that $X^\lambda$ is nonequidistant for all $\lambda>0$.  Then $\g(\la) \leq 1/(2\tau)$.
\end{thm}
Note that we do not require any stabilization or connectedness and also impose no direct restrictions on $\mu$ and on $\ell$. The proof of Theorem~\ref{thm-2degree1cluster} rests on showing absence of percolation in the SINR graph with a maximal degree given by $2$. The fact that SINR graphs with $\gamma>0$ have degrees less than $1+1/(\tau\gamma)$ is already stated in~\cite[Theorem 1]{DBT05}; an immediate consequence of this assertion is that there is no percolation in the case $\gamma \geq 1/\tau$ when degrees are at most 1. These claims can easily be seen to hold for any simple point process in any dimension, although in \cite{DBT05} only the case of a two-dimensional homogeneous PPP was considered. Theorem~\ref{thm-2degree1cluster} is the first improvement of this bound since then, applicable to stationary CPPs and thus in particular also covering the case of homogeneous PPPs in all dimensions.

\medskip
Finally, our third main result states that the critical intensity parameter for the SINR graph can be represented as the critical threshold for percolation of an associated Gilbert graph in any dimension. For this we assume a simpler setting in which $\L(\d x)$ equals the Lebesgue measure $\d x$, i.e., the CPP is in fact a PPP, and the powers are non-random and given by $p>0$.
% The associated SINR graph is denoted by $g_{(\g,P)}(X^\la)$ and correspondingly we write $\la_P$ for the critical intensity. Then, n
Note that for $\g=0$, the SINR graph is in fact a \emph{Poisson--Gilbert graph} (cf.~\cite{G61}) with connection radius given by 
\begin{align}\label{GilbThresh}
\rG=\ell^{-1}(\tau \No/p),
\end{align}
which is a well-defined quantity if $\ell(0)>\tau N_o/p$ and the conditions of Theorem~\ref{thm-randomfadings} on $\ell$ hold.

%We denote this Gilbert graph by $g_{\rG}(X^\lambda)$. Here, we r
Recall that the Gilbert graph based on a simple point process $Y$ with connection radius $r>0$ has vertex set $Y$ and an edge between two different points of $Y$ whenever the distance of the two points is less than $r$, and the name Poisson--Gilbert graph corresponds to the case when $Y$ is a homogeneous PPP. It is a standard result in continuum percolation that for the Poisson--Gilbert graph with connection radius $0<r<\infty$ in $d \geq 2$ dimensions, there exists a unique critical intensity $0<\lc(r)<\infty$ that separates a supercritical regime, where $\la>\lc(r)$, in which the Gilbert graph percolates with probability one and a subcritical regime, where $\la<\lc(r)$, in which the Gilbert graph does not percolate almost surely, cf.~for example~\cite[Section 3]{MR96}.

\begin{thm}\label{thm-PPPlambdaequality}
Let $d \geq 2$, $\No,\tau,p>0$, $\L(\d x)=\d x$ and $\ell$ be well-behaved with $\ell(0) > \tau \No/p$, then $\lacc=\lc(\rG)$.
\end{thm}
Theorem~\ref{thm-PPPlambdaequality} extends the result~\cite[Theorem 1]{DF06} to dimensions $d\ge3$ using new techniques, see Section~\ref{sec-StratProof} for details.

In the following section we present our main technical conditions together with examples for which our main theorems are applicable.

%===========================================
%===========================================
%===========================================

\section{Stabilization, asymptotic essential connectedness and examples}\label{sec-Ex}
The following definitions were recently introduced in~\cite{HJC17} in order to prove existence of a unique non-trivial critical intensity threshold for Cox--Gilbert graphs with fixed connectivity threshold. Let us recall that $Q_n=[-n/2,n/2]^d$ for $n>0$ and $d \in \N$, let $Q_n(x)=Q_n+x$ denote the box with side length $n$, centered at $x\in\R^d$, and $\dist(x,A):=\inf\{|x-y|\colon y\in A\}$ for $x \in \R^d$ and $A \subset \R^d$. We start with the definition of stabilization that can be understood as a quantitative spatial mixing property of the directing measure of a CPP. 
\begin{defn}[Stabilization]\label{defn-stabilization}
The random measure $\Lambda$ is called \emph{stabilizing} if there exists a random field of \emph{stabilization radii} $R=\{ R_x \}_{x\in \R^d}$ defined on the same probability space as $\Lambda$ such that, writing
\[ R(Q_n(x)) = \sup_{y \in Q_n(x) \cap\,  \Q^d} R_{y} ,\qquad n \geq 1, ~x \in \R^d, \]
the following hold.
\begin{enumerate}
\item\label{first-stabilizing} $(\Lambda,R)$ is jointly stationary,
\item $\lim_{n \uparrow \infty} \mathbb P(R(Q_n)<n)=1$, and
\item for all $n \geq 1$, non-negative bounded measurable functions $f$, and finite $\varphi \subset \R^d$ with $\dist(x,\varphi\setminus \{ x \})>3n$ for all $x \in \varphi$, the following random variables are independent:
\[  f(\Lambda_{Q_n(x)}) \mathds 1 \{ R(Q_n(x)) <n \}, \qquad x \in \varphi, \] 
\end{enumerate}
where for a measurable set $A \subseteq \R^d$, $\Lambda_A$ denotes the restriction of $\Lambda$ to $A$.
\end{defn}
A stronger form of stabilization is when $\Lambda$ is \emph{$b$-dependent}. That is, the restrictions $\Lambda_A$ and $\Lambda_B$ of $\Lambda$ to the measurable sets $A,B\subset\R^d$ are independent whenever $\dist(A,B)>b$ for some $b>0$. For $b$-dependence of subsets of $\Z^d$ we will use the analogous definition but with $\dist$ replaced by the $\ell^{\infty}$-distance. Next we give a definition of asymptotic essential connectedness, a suitable way of capturing connectedness of the support of the directing measure of a CPP with high probability. 
\begin{defn}[Asymptotic essential connectedness]\label{defn-asessconn}
The stabilizing random measure $\Lambda$ with stabilization radii $R$ is \emph{asymptotically essentially connected} if for all $n \geq 1$, whenever $R(Q_{2n})<n/2$, we have that
\begin{enumerate}
    \item $\supp(\Lambda_{Q_n})$ contains a connected component of diameter at least $n/3$,
    \item any two connected components of  $\supp(\Lambda_{Q_n})$ of diameter at least $n/9$ are contained in the same connected component of $\supp(\Lambda_{Q_{2n}})$. 
\end{enumerate} 
\end{defn}

The class of stabilizing random measures includes a number of interesting and relevant examples, 
for instance directing measures given via random tessellations based on PPPs. As already proven in~\cite[Section 3.1]{HJC17}, for example the edge-length measures of {\em Poisson--Voronoi tessellations} are asymptotically essentially connected (and hence also stabilizing), and it was also pointed out there that the same property for {\em Poisson--Delaunay tessellations} can be proven very similarly. It is nevertheless easy to see that these intensity measures are not $b$-dependent for any $b>0$. However, let us note that the edge-length measures of {\em Poisson line tessellations} in $\R^2$ are not even stabilizing. 

Stabilizing random measures that are absolutely continuous with respect to the Lebesgue measure include the directing measure of some {\em modulated PPPs} or {\em shot-noise fields} with compactly supported kernel. For the purpose of the present paper, a modulated PPP is defined
with directing measure $\Lambda(\d x)=\lambda\one\{x\in \Xi\}\d x+ \lambda'\one\{x\notin \Xi\}\d x$, for some Poisson--Boolean model $\Xi$ with constant connection radii, where the definition of a Poisson-Boolean model (with constant connection radii) will be presented at the beginning of Section~\ref{sec-lambda=lambdaproof}, and $\lambda,\lambda'\ge0$. As noted in \cite[Section 2.1]{HJC17}, the intensity measure that this definition yields is easily seen to be $b$-dependent for some $b>0$, and if $\lambda$ and $\lambda'$ are positive, then $\Lambda$ is asymptotically essentially connected. There exist examples, both for $\lambda>0$ and $\lambda'=0$ as well as $\lambda=0$ and $\lambda'>0$, such that asymptotic essential connectedness fails, see \cite[Section 2.5.1]{T18} for details. However, if $\Xi$  is in the supercritical regime for percolation and $\lambda>0$, then $\Lambda$ is asymptotically essentially connected, which follows from \cite[Theorems 2 and 5]{PP96}, as was observed in \cite[Section 2.1]{HJC17}. The general definition of a modulated PPP can be found in \cite[Section 5.2.2]{CSK+13}; here, the construction is similar to the case presented here, but $\Xi$ can be a general random closed subset of $\R^d$, and hence the arising directing measure need not even be stabilizing, as explained in \cite[Section 2.5.1]{T18}. Finally, without going into details, let us mention that the case when $\Xi$ is a Poisson--Boolean model with random radii, it is possible that the corresponding directing measure is stabilizing but not $b$-dependent for any $b>0$, see \cite[Example 3.4]{JTC20}. 

{\em Shot-noise fields} have directing measures of the form $\Lambda(\d x)=\sum_{i\in I}\kappa(Y_i-x)\d x $, with $(Y_i)_{i\in I}$ a homogeneous PPP and $\kappa\colon \R^d\to[0,\infty)$ compactly supported, cf.~\cite[Example 2.2]{HJC17}. They are always $b$-dependent for some $b>0$ but not asymptotically essentially connected in general, see~\cite[Section 2.1]{HJC17}, but in some relevant cases they are, see~\cite[Section 2.5.1]{T18}.

\section{Methods and discussion}\label{sec-StratProof} 
In this section, we lay out the strategies for the proofs of our main results, and comment on limitations and further extensions of the statements presented.

\subsection{Strategy of proof and discussion for Theorem~\ref{thm-randomfadings}}\label{sec-StratProof1}
As mentioned in the introduction, the statement of Theorem~\ref{thm-randomfadings} is an extension of the results of~\cite{L19} to the case of stabilizing CPPs. For the proof, we combine the approach used for~\cite[Theorem 4.5]{L19} for handling random radii and the approach used for~\cite[Theorem 2.4]{T18} for dealing with the spatial correlations of the directing measure $\L$ of the CPP.  To begin with, an easy coupling argument, see~\cite[Section 4.2.3.4]{T19b}, implies that as long as the powers are bounded, all positive results of~\cite{T18} about percolation in the Cox-SINR graph for asymptotically essentially connected $\L$ are applicable. More precisely, we have the following proposition for the Cox-SINR graph with random bounded powers.  
%, where we recall that we always assume that $\P(P_0>0)>0$.
\begin{prop}{\cite[Theorem 2.4 and Proposition 2.7]{T18}}\label{prop-T18}
Let $d\ge 2$, $\No,\tau>0$, $\P(P_o>0)>0$, $\Lambda$ be stabilizing and $\ell$ be well-behaved.
If $\Psup<\infty$ and $\ell(0)>\tau\No/\Psup$, then $\lacc<\infty$ holds if at least one of the following conditions is satisfied: 
\begin{enumerate} 
\item\label{second-boundedfading} $\ell$ has unbounded support, $\Lambda$ is $b$-dependent for some $b>0$ and $\E[\exp(\alpha\Lambda(Q_1))]<\infty$ holds for some $\alpha>0$, and at least one of the following conditions hold: $\Lambda$ is asymptotically essentially connected, or $P_{\sup}$ is sufficiently large, or
\item\label{stab} $\ell$ has bounded support, and $\Lambda$ is asymptotically essentially connected, or 
\item\label{bddsupportstab} $\ell$ has bounded support, and $\sup \supp(\ell)$ and $P_{\sup}$ are both sufficiently large.
\end{enumerate}
\end{prop}
Note that we have formulated the Condition~\eqref{second-boundedfading} in Proposition~\ref{prop-T18} more generally than what was stated in~\cite{T18}, and the Condition~\eqref{bddsupportstab} does not appear in \cite{T18}. However, the proof from~\cite{T18} can also be adapted to these more general cases. Indeed, let us first explain how the case of a constant power $\mu=\delta_p$, $p>0$, can be handled using the methods of \cite{T18}. The cases where $\Lambda$ is asymptotically essentially connected in Proposition~\ref{prop-T18} for constant powers are covered by \cite[Theorem 2.4, Part (2)]{T18}. Further, the methods of the proof of \cite[Proposition 2.7]{T18} apply to the case when $\Lambda$ is stabilizing but not necessary asymptotically essentially connected. To see this, note that the arguments of that proof require the connection radii $r_{\rm B}$ (see \eqref{GilbThresh}) to be large enough. Now, if $\supp(\ell)$ is unbounded, then one can always make $r_{\rm B}$ arbitrarily large via choosing the power value $p$ sufficiently large, which corresponds to the case of large $P_{\sup}=p$ in \eqref{second-boundedfading}. Else, this is not always possible because $\sup_{p>0} \ell^{-1}(\tau N_0/p)$ equals the finite number $\sup \supp(\ell)$. However, once $\sup \supp(\ell)$ is sufficiently large, one can make $r_{\rm B}$ sufficiently large such that the proof of \cite[Proposition 2.7]{T18} becomes applicable. Hence, we see that Proposition~\ref{prop-T18} indeed follows from \cite{T18} for fixed $p>0$. Now, if $\mu$ is not concentrated at one point, then one can always choose $p_2 \geq p_1>0$ such that $\mu([p_1,\infty))>0$ and $\mu((p_2,\infty))=0$. Then, if $p_1$ is sufficiently large, then the above arguments imply that there exists an infinite connected component in the subgraph of the SINR graph spanned by all vertices $X_i$ where $i \in I$ is such that $P_i>p_1$, for all sufficiently large $\lambda>0$ and all sufficiently small $\gamma>0$. Here, we bound all power values corresponding to the interferences by $p_2$ from above. See Section~\ref{sec-proof_thm1b}, Step~\ref{step-subgraph} for further details of a very similar argument.
%Indeed, if the SINR graph with constant powers $p_1$ percolates for some intensity $\lambda$ and interference-cancellation factor $\gamma>0$, then the SINR graph with random radii percolates if the intensity is at least $\lambda/\mu([p_1,\infty))$ and the interference-cancellation factor is at most $\gamma p_1/p_2$. Indeed, for a transmission between two points $X_i,X_j \in X^\lambda$ with $P_i,P_j>p_1$, the transmitted power is at least $p_1$ and the powers appearing in the interference are at most $p_2/p_1$ times this quantity. Now, the process of transmitters with powers at least $p_1$ is an independent thinning of $X^\lambda$ with survival probability $\mu([p_1,\infty))$, and hence itself a Cox point process with directing measure $\lambda\Lambda$, cf~\cite[Colouring Theorem and Marking Theorem]{K93}. The aforementioned choice of parameters guarantees that already the subgraph of the SINR graph spanned by such transmitters percolates, and hence so does the entire SINR graph. 
We conclude that $\lacc<\infty$ holds under the assumption of Proposition~\ref{prop-T18}. Given Proposition~\ref{prop-T18}, in the present paper it suffices to prove the case when $\Psup=\infty$. We prove Theorem~\ref{thm-randomfadings} in Section~\ref{sec-Proof1}. 

%Bene does not understand: Since for $\Lambda$ stabilizing, absence of a supercritical phase in the Cox--Gilbert graph only occurs for small $r$, which corresponds to small $\esssup~\zeta$ and sufficiently large $\sup~\supp~\ell$ in the SINR setting, it will be sufficient to require stabilization of $\Lambda$ in our assertions. 

\medskip
Let us comment on some further aspects of Theorem~\ref{thm-randomfadings}. First, as for Condition \eqref{first-fading} in Theorem~\ref{thm-randomfadings}, an extension to the general stabilizing case is not possible in general. Indeed, even if $P_o$ has very heavy tails, as soon as $\supp(\ell)$ is bounded, the radii of the associated Cox--Gilbert graph with random radii are bounded. Then, it is not hard to exhibit examples of stabilizing directing measures $\L$, such that $\lacc=\infty$, see the examples in \cite[Section 2.5.1]{T18}. 

Second, if $\Lambda$ is such that $\Lambda(Q_1)$ is almost surely bounded, then the exponential-moment condition  \[ \E[\exp(\alpha\Lambda(Q_1))]<\infty \numberthis\label{expmomentscond} \] of condition \eqref{second-fading} in Theorem~\ref{thm-randomfadings} clearly holds for all $\alpha >0$. E.g., this is the case for the modulated PPP with $\lambda,\lambda' \geq 0$. Further, \eqref{expmomentscond} holds for shot-noise fields for all $\alpha>0$, see e.g.~\cite[Section 2.5.1]{T18}. For Poisson--Voronoi and Poisson--Delaunay tessellations, the $b$-dependence assumption in \eqref{second-boundedfading} fails for all $b>0$, and hence percolation in the SINR graph can only be concluded for compactly supported $\ell$. On the other hand, it was verified in \cite{JT19} that for these two kinds of tessellations in two dimensions, $\E[\exp(\alpha\Lambda(Q_1))]<\infty$ holds for all $\alpha>0$; it is not known whether the same holds in higher dimensions. 

Third, the moment conditions on $P_o$ may look surprising at first. Indeed, why do we need to upper bound moments of $P_o$ in order to guarantee percolation in an SINR graph? This is indeed counterintuitive in view of the Gilbert graph since there larger radii would lead to better connectivity. However, in the SINR graph, as mentioned above, larger powers also increase interference and thus also might decrease connectivity. The classical approach used in~\cite{DF06, BY13, T18} to establish percolation in SINR graphs is to show that the underlying Gilbert graph satisfies some strong connectivity properties and at the same time the interferences can be uniformly bounded on large connected areas. We follow this approach as well, however, the random powers dictate several workarounds.
%, we also give more details in Section~\ref{sec-proof_comment}.  

Fourth, the Condition~\eqref{second-fading} in Theorem~\ref{thm-randomfadings} is not necessarily optimal. However, we believe that if percolation with unbounded $\supp(\ell)$ and without exponential moments of $P_o$ is possible, then the proof for this statement must be rather different from ours. An interference-control argument may not be possible at all, instead one should be able to show that the SINR values are sufficiently large for many transitions yielding satisfactory connectivity of the network for percolation. Let us mention a similar problem. It was conjectured in \cite{DBT05} that in the case with constant powers, in order to have percolation in the SINR graph for large $\la$, $\ell$ has only to have integrable tails but it may explode at zero. However, the setting where $\lim_{r \downarrow 0} \ell(r)=\infty$ is such that the classical interference-control argument, as exhibited in~\cite{DF06} and adapted to the case of random powers in Section~\ref{sec-proof_thm1b}, certainly cannot work. Indeed, the interferences are almost-surely finite but they have infinite expectation, see~\cite{D71}, hence there is no hope to apply a version of the exponential Markov inequality. Let us also note that the results of~\cite{D71} also imply that, if the tails of $\ell$ are not integrable, then SINR graphs with $\gamma>0$ have no edges. 

Finally, under the assumptions of Theorem~\ref{thm-randomfadings} on $\ell$, for $\gamma=0$, the SINR graph $G_{0}(\mathbf X^\lambda)$ is a Gilbert graph with i.i.d.\ random radii $R_i=\ell^{-1}(\tau N_o/P_i)$. Let $R_o$ denote a generic random variable having the same distribution as $R_i$. Now, if all other parameters are kept fixed, it is easy to see that $\gamma \mapsto \P(\Gg \text{ percolates})$ is decreasing. Hence, if almost surely there is no percolation in the SINR graph for $\gamma=0$, then the same holds for all $\gamma>0$. Further, as already mentioned, percolation properties of Gilbert graphs can equivalently be expressed in terms of the corresponding Boolean models. This way, the recent result \cite[Theorem 2.6, Part (2)]{JTC20} about existence of a subcritical phase in Cox--Boolean models immediately implies the following assertion. If $\Lambda$ is $\phi$-stabilizing and $R_o$ is unbounded with $\E[R_o^d]<\infty$, then $\lacc>0$. Here, the notion of $\phi$-stabilization (cf.~\cite[Definition 2.5]{JTC20}) is very similar to our definition of stabilization, and many relevant stabilizing examples are also $\phi$-stabilizing. This observation complements the result $\lacc<\infty$ in Theorem~\ref{thm-randomfadings}. Moreover, it improves the assertion of \cite[Section 4.2.3.4]{T19b} that $\lacc<\infty$ holds for bounded $P_o$ (equivalently, bounded $R_o$) if $\Lambda$ is stabilizing and $\ell$ satisfies the assumptions of Theorem~\ref{thm-randomfadings} with $\ell(0) > \tau \No/\essinf \mu$ or $\ell(0) \leq \tau \No/\Psup$. 

%=================================================
%=================================================

\subsection{Strategy of proof and discussion for Theorem~\ref{thm-2degree1cluster}}\label{sec-strategy2}\label{sec-StratProof2}
\begin{figure}[!htpb]
\centering
\input{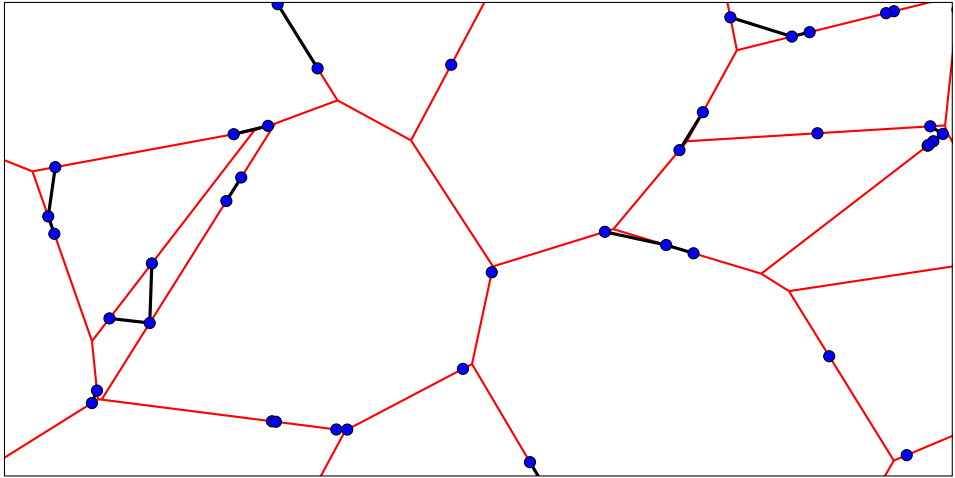}
\caption{A typical realization of a Cox-SINR graph (with blue vertices and black edges) with directing measure given by the edge-length measure of a two-dimensional Poisson--Voronoi tessellation (in red) in a box, with $\No=\tau=1$, constant powers equal to $1$, and a suitable path-loss function $\ell$. The interference-cancellation factor is set to $\gamma=1/(2\tau)$. We see only a few vertices having degree two, the largest connected component is of size three, and there are no cycles in the graph. As indicated by Proposition~\ref{prop-2degree1cluster} the graph is highly disconnected. 
%For the definition of the Poisson--Voronoi tessellation and its applications in modelling urban street systems, see~\cite{HJC17} and the references therein.
}\label{figure-2degree}
\end{figure}

We call a maximal connected component in a graph a \emph{cluster}. As already pointed out in~\cite[Theorem 1]{DBT05}, for $\gamma>0$, all degrees in $\Gg$, where $X^\la$ is a PPP, are less than $1+1/(\tau\gamma)$ for any choice of $\la,\tau>0$ and $\No \geq 0$. In other words, all vertices in $\Gg$ have at most $1+1/(\tau\gamma)$ neighbors. It is not hard to see that this property remains true if the PPP is replaced by a CPP, or even any simple point process, see~\cite[Section A.3]{T19b}.
% further, for any choice of the nonnegative random power distribution $\zeta=\P \circ P_0^{-1}$. This holds for a general simple point process $X^\lambda$ and a general nonnegative $\zeta$, more precisely, for given $X_i \in X^\lambda$, the number of indices $j \neq i$ such that $\SINR(X_j,X_i,\mathbb X^\lambda)>\tau$ is less than $1+\frac{1}{\tau\gamma}$, see \cite[Section A.3]{T19b}. All these assertions depend on the combinatorial structure of the notion~\eqref{SINRfading} of SINR and are entirely non-spatial. (Note that if $\SINR(X_j,X_i,\mathbb X^\lambda)>\tau$, then $P_j>0$, even for $N_0=0$.) \\%Thus, if $\gamma$ is relatively large compared to $1/\tau$, then all degrees in the SINR graph are very low, and therefore it is less straightforward to show that strong local connectivity of the graph can be guaranteed with positive probability than for $\gamma>0$ very close to zero. In particular, we need at least partially new methods in order to verify that there is at most one infinite cluster in the graph.
Thanks to the degree bounds, any such Cox-SINR graph with random powers for which $\gamma \geq 1/\tau$ has no infinite cluster since it has degrees bounded by 1. For $\gamma \in [1/(2\tau),1/\tau)$, we have an \emph{a priori} degree bound of 2, which implies that all maximal connected components of SINR graphs are finite cycles or paths that are infinite in zero, one or two directions. This reminds of a one-dimensional percolation model, and thus the conjecture is that it contains no infinite clusters under general assumptions on the directing measure of the CPP, see Figure~\ref{figure-2degree} for an illustration. 
The following proposition shows that this is indeed true for the Cox-SINR graph with random powers. 
\begin{prop}\label{prop-2degree1cluster}
Let $d\ge 1$, $\No\geq 0$, $\tau>0$ and $\gamma \geq 1/(2\tau)$, then for $\Lambda$ nonequidistant, 
\[ \P(\Gg\text{ percolates})=0. \]
\end{prop}
The statement of Theorem~\ref{thm-2degree1cluster} is an immediate consequence of Proposition~\ref{prop-2degree1cluster}, the proof of which can be found in Section~\ref{sec-twodegreeproof}. The proof employs a fine configuration-wise analysis of the SINR graph, which seems to be new in the literature. Moreover, we expect the proof to hold for SINR graphs based on a large class of simple nonequidistant stationary point processes.

Let us comment on a further aspect of Theorem~\ref{thm-2degree1cluster}. It can be observed that the proof of Proposition~\ref{prop-2degree1cluster} does not use the precise numerical relation $\gamma \geq 1/(2\tau)$, but rather just the fact that the SINR graph has degrees at most 2. Hence, once $\Lambda$ is nonequidistant, the result holds as soon as the SINR graph has degrees bounded by 2. Note that for example if $N_o>0$, $P_o$ bounded and $\ell$ continuous on $[0,\infty)$, then one can derive a stricter upper bound on the degrees (depending on several parameters) along the lines of the proof of \cite[Theorem 1]{DBT05}.
%Its proof can be found in Section~\ref{sec-twodegreeproof}. \color{red} Apart from checking the proof, the main question is whether we can do this for random radii. See my remarks after the proof. \color{black}
%
%
%
% Let us note that for two-dimensional Poisson point processes we have
% \[ \lambda_{N_0,\tau,P} = \inf \{ \lambda>0 \colon \gamma^*(\lambda)>0 \}, \numberthis\label{lambdacrittwodefinitions}\]
%{\red reference missing,} and it is plausible to conjecture that the same holds for a large class of point processes. Note that for any ergodic point process, the probability of percolation of the SINR graph is in $\{ 0,1 \}$ (also for $\gamma=0$), cf.~\cite[Section 2]{MR96}. 
% 

%=================================================
%=================================================

\subsection{Strategy of proof for Theorem~\ref{thm-PPPlambdaequality}}\label{sec-StratProof3}
As mentioned previously, we have $G_{0}(X^\lambda)=g_{\rG}(X^\lambda)$ for all $\lambda>0$ in the Poisson-SINR graph with fixed powers $p$, where $r_{\rm B}$ is defined in~\eqref{GilbThresh}. We use $X^\la$ instead of $\X^\la$ since the marks are non-random. 
Moreover, note that the increase of the interference-cancellation factor $\g$ can only lead to edges being removed from the graph and hence 
there is a monotonicity of $G_{\g}(X^\lambda)$ with respect to $\g$. Additionally, there is a monotonicity of $g_{\rG}(X^\lambda)$ with respect to $\la$, which together implies that $\lacc\geq \lc(\rG)$. We have the following equivalence result from~\cite{DF06} for the two-dimensional Poisson-SINR graphs. 
\begin{thm}{\cite{DF06}}\label{thm-d=2lambda=lambda}
Let $d = 2$, $\No,\tau,p>0$, $\L(\d x)=\d x$, and $\ell$ be well-behaved. Then $\lacc=\lc(\rG)$.
\end{thm}
In words, this result states that for any $\lambda>0$ such that the Poisson--Gilbert graph $g_{\rG}(X^\lambda)$ is supercritical, there exists $\gamma>0$ such that also the Poisson-SINR graph $G_{\g}(X^\lambda)$ percolates. In an extended context of SINR graphs, it was shown that this percolation is preserved if the transmitters forming a PPP experience additional interference coming from a weakly $\alpha$-sub-PPP, see~\cite[Section 3.4]{BY13}. 

The proof of Theorem~\ref{thm-d=2lambda=lambda} employs Russo--Seymour--Welsh type arguments for the Poisson--Gilbert graph in two dimensions, see~\cite[Section 4]{MR96} and~\cite[Section 3]{DF06}. These arguments do not have a known analogue in the Poisson case for $d \geq 3$, or in the general Cox case even for $d=2$. Note that the results of \cite{T18} only imply that $\lacc<\infty$ for $d \geq 3$ and $\L(\d x)=\d x$. However,~\cite[Section 2.1]{HJC17} includes some further observations about Gilbert graphs in $d \geq 3$ dimensions, originating from results of~\cite{PP96}, that allow us to conclude the analogue of Theorem~\ref{thm-d=2lambda=lambda} for the higher-dimensional Poisson case. The proof of Theorem~\ref{thm-PPPlambdaequality} will be carried out in Section~\ref{sec-lambda=lambdaproof}, where we will also recall the results of~\cite{PP96} that are relevant for our analysis.

%===========================================
%===========================================
%===========================================

\section{Proofs}\label{sec-Proofs}
For the proofs it will be convenient to define the SINR of $X_i,X_j \in X^\la$, $X_i \neq X_j$, via
\[ \SINR(X_i,X_j,\X^\la) = \frac{P_i\ell(|X_i-X_j|)}{\No+\gamma \sum_{k \in I \setminus \{ i, j\}} P_k \ell(|X_k-X_j|)}. \numberthis\label{SINRproof}\]
%=================================================
%=================================================

\subsection{Proof of Theorem~\ref{thm-randomfadings}}\label{sec-phasetransitionproof}\label{sec-Proof1}
Let us first carry out the proof under Condition~\eqref{second-fading} of the theorem in Section~\ref{sec-proof_thm1b}. The proof under the Conditions~\eqref{first-fading} and~\eqref{third-fading} of the theorem are presented in Sections~\ref{sec-proof_thm1a} and~\ref{sec-proof_thm1c}, respectively. 
%=================================================

\subsubsection{Proof of Theorem~\ref{thm-randomfadings} Part~\eqref{second-fading}}\label{sec-proof_thm1b}
%As discussed in Section~\ref{sec-StratProof1}, it suffices to consider the case $\Psup=\infty$. 
For fixed $\lambda$ and $\gamma$, in order to show that $\Gg$ percolates, it suffices to verify that a subgraph of it contains an infinite cluster. Our proof consists of four steps. First, for $\gamma,\lambda>0$, we define a subgraph that is included in a Cox--Gilbert graph (with constant connection radii). Second, we map this subgraph to a lattice percolation model and show that this discrete model percolates for large $\lambda$ for a suitable choice of auxiliary parameters. In particular, since $\Lambda$ is only assumed stabilizing, the connection radius of the Gilbert graph must be large enough so that the graph percolates for large $\lambda$. In this step, we are able to employ multiple arguments of \cite{DF06,HJC17,T18}. Our interference-control assertion, Proposition~\ref{prop-interferencecontrol}, is presented here. Third, using the subgraph, we make a choice of $\gamma>0$ such that percolation in the discrete model implies percolation in the SINR graph $\Gg$, which is done analogously to \cite{DF06}. Fourth, we carry out the proof of Proposition~\ref{prop-interferencecontrol}, combining arguments of \cite{DF06,T18} for SINR graphs with constant powers and arguments used in~\cite{L19} for Poisson-SINR graphs with random powers.

%=================================================

\begin{step}\label{step-subgraph}
A subgraph of the SINR graph.
\end{step}
We first present a general construction of a subgraph of $\Gg$ for $\gamma,\lambda>0$. Let $r_o>d_o$. Since both $P_o$ and $\supp(\ell)$ are unbounded, we have
\[ p(r_o): = \P\big( \ell^{-1}\big(\tau \No/P_o\big) \geq r_o \big) = \P\Big( P_o \geq\tau \No/\ell(r_o) \Big) >0. \]
Let us define the independent thinning
\[ X^{\lambda,-} := \{ X_i \in X^\lambda \colon P_o \geq \tau \No/\ell(r_o) \} \]
of $X^\lambda$ with survival probability $p(r_o)$. According to \cite[Colouring Theorem]{K93}, $X^{\lambda,-}$ is a CPP with directing measure $\lambda p(r_o)\Lambda$. Now, let us define a subgraph $G^-_{\g}(\X^\la)$ of $\Gg$ as follows. The vertex set is $X^{\lambda,-}$, and two vertices $X_i,X_j \in X^{\lambda,-}$, $X_i \neq X_j$, are connected by an edge if and only if
\[ \SINR^{-}(X_i,X_j,\X^\la): = \frac{\big(\tau \No/\ell(r_o)\big) \ell(|X_i-X_j|)}{\No+\gamma \sum_{k \in I \setminus \{ i, j\}} P_k \ell(|X_k-X_j|)} > \tau \numberthis\label{SINRminus} \]
and the analogously defined $\SINR^-(X_j,X_i,\X^\la)$ also exceeds $\tau$.
Note that for $X_i,X_j \in X^{\lambda,-}$, in the numerator of $\SINR(X_i,X_j,\X^\la)$, for the power of $X_i$ we have $P_i \geq \tau \No/\ell(r_o)$, whereas the denominators of \eqref{SINRproof} and \eqref{SINRminus} are equal, and the same holds with the roles of $i$ and $j$ interchanged. Hence, $G^-_{\g}(\X^\la)$ is indeed a subgraph of $\Gg$ for any $\gamma \geq 0$. As for $\gamma=0$, $G^-_{0}(\X^\la)$ equals the Cox--Gilbert graph $g_{r_o}(X^{\lambda,-})$ with connection radius $r_o$ and vertex set $X^{\lambda,-}$. In words, in order to obtain  $G^-_{\g}(\X^\la)$ from $\Gg$, one first thins out vertices with small powers, in order to get rid of vertices with small values of the connection radius $\rG^i$, where
\begin{align}\label{r_B}
\rG^i=\ell^{-1}(\tau \No/ P_i).
\end{align}
Then, one bounds the powers of the remaining vertices by $\tau \No/\ell(r_o)$ from below.

\begin{step}\label{step-lattice}
Mapping the subgraph to a lattice-percolation problem and percolation on the lattice.
\end{step}

Now we are in a position to adapt to the setting of \cite[Section 3.2.2]{T18} and use strong connectivity of $g_{r_o}(X^{\lambda,-})$ in case $r_o$ is sufficiently large and $\lambda$ is chosen according to $r_o$. 
Together with an interference-control argument presented below, this will allow us to verify Theorem~\ref{thm-randomfadings} Part~\eqref{second-fading}.

First, let us recall the definition of rescalings of a Gilbert graph that were also used in \cite{T18}. For $c>0$ and a Gilbert graph $G$ with connection radius $r>0$, deterministic vertex set $V \subset \R^d$, $d \geq 1$, and edge set $E=\{ (x,y) \in V\times V \colon x \neq y, |x-y|<r \}$, the graph $c G$ is defined with vertex set $cV=\{cx\colon x\in V\}$ and edge set $cE=\{ (cx,cy) \colon (x,y) \in E \}$. It is easy to see that $c G$ is a Gilbert graph with vertex set $cV$ and connection radius $cr$. For Gilbert graphs with random vertex sets (e.g., if the vertex set is given by a random simple point process), rescalings of the graph are defined realizationwise. From~\cite[Section 7.1, Proof of Theorem 2.9 (Convergence in Bounded Domains)]{HJC17} we know that in the coupled limit $\widetilde r \uparrow \infty, \widetilde \lambda \downarrow 0$ and $\widetilde \lambda \widetilde r^d = \widetilde \varrho>0$, we have that $\widetilde r^{-1} g_{\widetilde r}(X^{\widetilde \lambda})$ converges weakly to the graph $g_{1}(Y^{\widetilde \varrho})$, where $Y^{\varrho}$ is a homogeneous PPP with intensity $\varrho$. Let us note that in~\cite[Section 7.1]{HJC17} this convergence is formulated equivalently for the Boolean model and the crucial point is that the convergence is guaranteed only in compact domains. 
%This convergence is to be understood in the sense that if we represent the graph $\widetilde r^{-1} g_{\widetilde r}(X^{\widetilde \lambda, \widetilde r})$ as the union of the closed straight line segments between the pairs of its vertices connected by its edges, and we intersect the arising set with an arbitrary compact set $K \subset \R^d$, then this intersection converges to the intersection of the analogously defined set corresponding to the limiting graph $g_{1}(Y^{\widetilde \varrho})$ with $K$. \color{red} In what topology? \color{black}

%For $\varrho>0$, let $Y^{\varrho}$ be a PPP with intensity measure (directing measure) $\varrho\Leb$. 
Let $\varrho_{\rm c}(1)$ be such that the Poisson--Gilbert graph $g_{1}(Y^{\varrho_{\rm c}(1)})$ is critical. Then, due to the scale invariance of Poisson--Gilbert graphs \cite[Section 2.2]{MR96}, for $\varrho>\varrho_{\rm c}(1)$, we can choose a smaller intensity $\varrho'<\varrho$ such that $g_{1}(Y^{\varrho'})$ is still supercritical.   Now, for $r>d_o$, we define $r_o(r)=r(\varrho/\varrho')^{1/d}$, $\lambda(r)=\varrho' r^{-d} (p(r_o(r)))^{-1}$ and $p(r)=\tau \No/\ell(r_o(r))$. Noting that $g_r(X^{\lambda(r),-})$ is a Cox--Gilbert graph with connection radius $r$ and stabilizing intensity $p(r_o(r)) \lambda(r)=\varrho' r^{-d}$, we have that $r^{-1}g_r(X^{\lambda(r),-})$ converges to the supercritical graph $g_1(Y^{\varrho'})$ on compact sets, as $r$ tends to infinity.
%Moreover, we have $\ell(r_o(r))=\tau \No/P(r)$. %(Note that $r_o(r)$ is called $r_{\rm B}(r)$ in \cite{T18}.)

Further, recalling that $R$ denotes the stabilization radii of $\Lambda$, we put $R(Q)=\sup_{x \in Q \cap \Q^d} R_x$ for any measurable set $Q \subseteq \R^d$. 

Using these notions, we construct a renormalized percolation process on $\Z^d$ as follows. For $n \geq 1$ and $r>d_o$, the site $z \in \Z^d$ is \emph{$(r,n)$-good} if 
\begin{enumerate}
    \item\label{first-nrgood} $R(Q_{6rn}(rnz))<rn/2$, 
    \item\label{second-nrgood} $X^{\lambda(r),-} \cap Q_{rn}(rnz) \neq \emptyset$, and
    \item\label{third-nrgood} every $X_i,X_j \in X^{\lambda(r),-} \cap Q_{3rn}(rnz)$ are connected by a path in $g_r(X^{\lambda(r),-}) \cap Q_{6rn}(rnz)$. 
\end{enumerate}
The site $z \in \Z^d$ is \emph{$(r,n)$-bad} if it is not $(r,n)$-good. Note that the process of $(r,n)$-good sites is 7-dependent thanks to the definition of stabilization. The following lemma has been verified in \cite[Section 3.2.2]{T18}. However, since in \cite{T18} it was not formulated as a lemma and two different proofs are presented for $d=2$ and $d \geq 3$, we provide a self-contained proof here for the reader's convenience.
\begin{lem}{\cite{T18}}\label{lem-T18_1}
Assume that the general conditions of Theorem~\ref{thm-randomfadings} plus the conditions in Part~\eqref{second-fading} hold. Then, for all sufficiently large $\lambda>0$ and for all $n \geq 1$ and $r >d_o$ with $rn$ sufficiently large, there exists $q_A=q_A(\lambda,rn) <1$ such that for any $N \in \N$ and pairwise distinct $z_1,\ldots,z_N \in\Z^d$,  
\[ \P(z_1,\ldots,z_N\text{ are all $(r,n)$-bad}) \leq q_A^N. \numberthis\label{qAestimate} \]
Further, for any $\eps>0$, one can choose $\lambda$ and $rn$ sufficiently large such that $q_A <\eps$.
\end{lem}
\begin{proof}
For $z \in \Z^d$, we write $J_{n,r}(z)$ for the event that $z$ satisfies \eqref{second-nrgood} and \eqref{third-nrgood} in the definition of $(n,r)$-goodness. Then, for any $n, r$ under consideration, the process of $(n,r)$-good sites is 7-dependent by the definition of stabilization. Further, we write $F_n(z)$ for the event that in the definition of $(n,1)$-goodness, the PPP $Y^{\varrho'}$ with intensity $\varrho'=\lambda(r)r^d$ satisfies \eqref{second-nrgood} with $X^{\lambda(r),-}$ replaced by $Y^{\varrho'}$ and \eqref{third-nrgood} with $g_r(X^{\lambda(r),-})$ replaced by $g_1(Y^{\varrho'})$ everywhere. The probability of $F_n(z)$ is independent of the choice of $z$ and tends to 1 as $n \to\infty$ thanks to the arguments of \cite[Section 5.2]{HJC17}, since the constant directing measure of the PPP $Y^{\varrho'}$ is certainly asymptotically essentially connected. Using a union bound and the well-known scale invariance of Poisson--Gilbert graphs, namely that for $\widetilde\lambda,\widetilde r>0$, $\widetilde r^{-1} g_{\widetilde r} (X^{\widetilde\lambda})$ equals $g_1(X^{\widetilde \lambda \widetilde r^d})$ in distribution, we conclude that for $z \in \Z^d$,
\[ \mathbb P(z\text{ is }(n,r)\text{-bad}) \leq \mathbb P \Big( R(Q_{6nr}(nrz)) \geq nr/2 \Big) +\mathbb P(F_n(z)^{\mathrm c})+|\mathbb P(F_{n}(z)^{\mathrm c})-\mathbb P(J_{n,r}(z)^{\mathrm c})|, \]
which can be made arbitrarily close to zero by first choosing $n$ large and then $r$ large according to $n$, due to the weak convergence of $r^{-1} g_r(X^\lambda)$ to $g_1(Y^{\varrho'})$ on $Q_{6n}(nz)$ as $r \to \infty$, $\lambda(r) \to 0$, $r^d \lambda(r) = {\varrho'}$.

Hence, applying \cite[Theorem~0.0]{LSS97}, for all sufficiently large $n$ and large enough $r$ chosen according to $n$, the 7-dependent process of $(n,r)$-good sites is stochastically dominated from below by a supercritical independent site percolation process. Moreover, the probability of a site of the independent site percolation process being closed can be made arbitrarily close to 0 via further increasing $nr$. This implies the lemma. 
\end{proof}

%Roughly speaking, the idea of the proof of this assertion in \cite[Section 3.2.2]{T18} is to compare $g_{r}(X^{\lambda,-})$ to a supercritical Poisson--Gilbert graph in the limit $\lambda \downarrow 0, r \uparrow \infty, \lambda r^d =\varrho'$, which is the weak limit of $g_{r_o}(X^{\lambda,-})$ if one considers restrictions of the graphs to compact sets \cite[Section 7.1]{HJC17}. 

We further proceed similarly to \cite{DF06,T18} by defining `shifted' versions of the path-loss function $\ell$. For $a\ge 0$, define \[ \ell_{a}(r) = \ell(0) \mathds 1 \big\{ r< a\sqrt d/2 \big\} + \ell\big( r-a\sqrt d/2  \big) \mathds 1 \big\{ r \geq  a\sqrt d/2 \big\}. \numberthis\label{elladef} \]
Note that $\ell_0=\ell$. 
Now, we define the\index{shot-noise process} shot-noise processes
\begin{eqnarray*}
& I_{a}(x)=\sum_{i \in I} P_i \ell_{a}(|x-X_i|), \qquad I(x)=\sum_{i \in I} P_i \ell(|x-X_i|), \qquad x \in \R^d,
\end{eqnarray*}
and note that $I_{0}(x)=I(x)$. By the triangle inequality, for $a \geq 0$, $I(x) \leq I_{a}(z)$ holds for any $z \in \R^d$ and $x \in Q_a(z)$. %Further, the following assertion can be verified analogously to \cite[Lemma 3.1]{T18}.
%\begin{lem}\label{lemma-ellshiftrandomradii}
%Let $m \geq 1$, $a \geq m$, $k \in \N$, $z \in k\Z^d$. Then for any $y \in Q_k(z)^{o} \cap \Z^d$ and $S \subseteq \R^d$ measurable, we have almost surely
%\[ \sum_{X_j \in X^\lambda \cap S} P_j \ell_{ka}(|mz-X_j|) \geq \sum_{X_j \in X^\lambda \cap S} P_j \ell_{a}(|m y-X_j|) 
%%I_{ka}(nz) \geq I_a(n z_i), \qquad \forall i \in \big[k^d\big]
%. \]
%\end{lem}
Now, the interference-control argument consists in verifying the following proposition. For $z \in \Z^d$, let us write $B_{r,n,M}(z)=\{ I_{6rn}(rnz) \leq M \}$.  
\begin{prop}\label{prop-interferencecontrol}
Assume that the general conditions of Theorem~\ref{thm-randomfadings} plus Part~\eqref{second-fading} hold. Then, for all $\lambda>0$, for all $n \geq 1$ and $r>d_o$ with $rn$ sufficiently large and for all $M>0$ sufficiently large, there exists $q_B=q_B(\lambda,rn,N) <1$ such that for all $N \in \N$ and for all pairwise distinct $z_1,\ldots,z_N \in \Z^d$ we have
\[ \P(B_{r,n,M}(z_1)^{\rm c} \cap \ldots \cap B_{r,n,M}(z_N)^{\rm c}) \leq q_B^N. \numberthis\label{qBstabilizing} \]
Further, for any $\eps>0$ and $\lambda>0$, one can choose $rn$ and $M$ sufficiently large such that $q_B < \eps$.
\end{prop}
The proof of this proposition is postponed until Step~\ref{step-interferencecontrol}. 
Once we have shown Proposition~\ref{prop-interferencecontrol}, we can derive the following corollary using a standard argument (see e.g.~the proof of \cite[Proposition 3]{DF06} or the one of \cite[Proposition 3.1]{T18}). For $z \in \Z^d$ let us define $C_{r,n,M}(z)=\{ z\text{ is }(r,n)\text{-good} \}\cap \{ I_{6rn}(rnz) \leq M \}$. 
\begin{cor}\label{prop-twoPeierls}
Assume that the general conditions of Theorem~\ref{thm-randomfadings} plus Part~\eqref{second-fading} hold. Then, for all sufficiently large $\lambda>0$, for all $r>d_o$ and $n \geq 1$ with $rn$ sufficiently large and for all $M>0$ sufficiently large, there exists $q_C =q_C(\lambda,rn,M)<1$ such that for all $N \in \N$ and for all pairwise distinct $z_1,\ldots,z_N \in \Z^d$ we have
\[ \P(C_{r,n,M}(z_1)^{\rm c} \cap \ldots \cap C_{r,n,M}(z_N)^{\rm c}) \leq q_C^N. \]
Further, for any $\eps>0$, one can choose $\lambda,rn,M$ sufficiently large such that $q_C < \eps$.
\end{cor}
\begin{step}\label{step-alsoSINRperc}
Percolation in the subgraph of the SINR graph.
\end{step}
Having Corollary~\ref{prop-twoPeierls} and employing a Peierls argument (cf.~\cite[Section~1.4]{G99}), we conclude that for $\lambda, rn,M$ sufficiently large, the process of $(r,n)$-good sites $z \in \Z^d$ such that $I_{6rn}(rnz) \leq M$ percolates. Thanks to the exact same arguments as in \cite[Section 5.2, Proof of Theorem 2.6]{HJC17}, this implies percolation of the Cox--Gilbert graph $G^-_{0}(\X^{\lambda(r)})=g_{r_o(r)}(X^{\lambda(r),-})$.  From this point of the proof it is classical to derive that $G^-_{\g}(\X^{\lambda(r)})$ percolates for small $\gamma>0$, see~\cite[Section 3.3]{DF06}. For the convenience of the reader, let us give the details here. We define 
\[ \gamma' = \frac{\No}{p(r)M} \Big( \frac{\ell(r)}{\ell(r_o(r))} - 1 \Big) =\frac{\ell(r_o(r))}{\tau M}\Big( \frac{\ell(r)}{\ell(r_o(r))} - 1 \Big) >0, \] 
where the strict inequality holds because $r_o(r)>r>d_o$ and $\ell$ has unbounded support. Then we have
\[ \frac{p(r) \ell(r)}{\No+\gamma' p(r) M} = \tau. \]
Now, let $X_i,X_j \in X^{\lambda(r),-}$ be situated in $Q_{rn}(rnz)$ respectively $Q_{rn}(rnz')$ for some sites $z,z' \in\Z^d$ included in the same infinite cluster of the process of $(r,n)$-good sites $z \in \Z^d$ satisfying $I_{6rn}(rnz) \leq M$ such that $|X_i-X_j| < r$. Then, for $\gamma<\gamma'$, we have
\[ \SINR(X_i,X_j,\mathbf X^\lambda) \geq \SINR^{-}(X_i,X_j,\mathbf X^\lambda) >  \frac{p(r) \ell(r)}{\No+\gamma' p(r) M} = \tau. \]
Thus, $X_i$ and $X_j$ are connected by an edge in $G^-_{\g}(\X^\la)$. Hence, $\Gg$ also percolates. Thus, we can conclude Theorem~\ref{thm-randomfadings} as soon as we have verified Proposition~\ref{prop-interferencecontrol}.
\begin{step}\label{step-interferencecontrol}
Proof of Proposition~\ref{prop-interferencecontrol}: the interference-control argument.
\end{step}

Similarly to \cite[Section 3.1.1]{T18}, we split the interference into two parts. For $x \in \R^d$, $n \geq 1$ and $r>0$, we put
\begin{align*}
& I_{6rn}^{\mathrm{in}} (x) = \sum_{X_i \in X^\lambda \cap Q_{12rn\sqrt d}(x)} \ell_{6rn}(|X_i-x|),~ I_{6rn}^{\mathrm{out}} (x) = \sum_{X_i \in X^\lambda \setminus Q_{12rn\sqrt d}(x)} \ell_{6rn}(|X_i-x|).
\end{align*}
Then, for $M>0$, if $I_{6rn}(x)>M$, then $I_{6rn}^{\mathrm{in}} (x)>M/2$ or $I_{6rn}^{\mathrm{out}} (x)>M/2$. Using a union bound and the fact that in Proposition~\ref{prop-interferencecontrol}, $M$ can be chosen arbitrarily large, it suffices to conclude the proposition both with $B_{r,n,M}(z_i)$ replaced by $B_{r,n,M}^{\mathrm{in}}(z_i)$ and with $B_{r,n,M}(z_i)$ replaced by $B_{r,n,M}^{\mathrm{out}}(z_i)$ everywhere in \eqref{qBstabilizing} for all $i \in \{1,\ldots,N\}$, where for $z \in \Z^d$ we write $B_{r,n,M}^{\mathrm{in}}(z)= \{ I_{6rn}^{\mathrm{in}}(rnz) \le M \}$ and $B_{r,n,M}^{\mathrm{out}}(z)= \{ I_{6rn}^{\mathrm{out}}(rnz) \le M \}$. Indeed, having these assertions, we can combine them similarly to Corollary~\ref{prop-twoPeierls}. 

We now verify Proposition~\ref{prop-interferencecontrol} with  $B_{r,n,M}(\cdot)$ replaced by $B_{r,n,M}^{\mathrm{in}}(\cdot)$ everywhere. For this assertion, instead of the assumption that $P_o$ and $\Lambda(Q_1)$ have some exponential moments, it suffices if they have a first moment (for $\Lambda(Q_1)$ this is automatic since $\E[\Lambda(Q_1)]=1$ by assumption). To be more precise, we prove the following lemma.
\begin{lem}\label{prop-interferencecontrolIN}
Assume that the general conditions of Theorem~\ref{thm-randomfadings} plus Part~\eqref{second-fading} hold. Further, let $\Lambda$ be stabilizing and $\E[P_o]<\infty$. Then, for all $\lambda>0$, for all $n \geq 1$ and $r>d_o$ with $rn$ sufficiently large and for all $M>0$ sufficiently large, there exists $q_B=q_B(\lambda,rn,N) <1$ such that for all $N \in \N$ and for all pairwise distinct $z_1,\ldots,z_N \in \Z^d$ we have
\[ \P(B_{r,n,M}^{\mathrm{in}}(z_1)^{\rm c} \cap \ldots \cap B_{r,n,M}^{\mathrm{in}}(z_N)^{\rm c}) \leq q_B^N. \numberthis\label{qBstabilizingIN} \]
Further, for any $\eps>0$ and $\lambda>0$, one can choose $rn$ and $M$ sufficiently large such that $q_B < \eps$.
\end{lem}
\begin{proof}
We use the following auxiliary discrete percolation process. A site $z \in \Z^d$ is \emph{$(r,n)$-tame} if
\begin{enumerate}
\item\label{first-rntame} $R(Q_{12rn\sqrt d}(rnz))<rn/2$, and
\item\label{second-rntame} $I_{6rn}^{\mathrm{in}}(rnz)\leq M $.
\end{enumerate}
A site $z \in \Z^d$ is \emph{$(r,n)$-wild} if it is not $(r,n)$-tame. The process of $(r,n)$-tame sites is $\lceil 12 \sqrt d + 1\rceil$-dependent according to the definition of stabilization. Thus, it follows from dependent-percolation theory \cite[Theorem~0.0]{LSS97} that, in order to verify Lemma~\ref{prop-interferencecontrolIN}, it suffices to show that for all $\lambda>0$, $\mathbb P(o\text{ is $(r,n)$-wild})$ can be made arbitrarily close to zero by choosing first $rn$ sufficiently large and then $M$ large enough accordingly. We have
\[ \mathbb P(o\text{ is $(r,n)$-wild}) \leq \mathbb P(R(Q_{12rn\sqrt d}(rnz)) \geq rn/2) + \mathbb P(I_{6rn}^{\mathrm{in}}(rnz)> M). \]
The first term can be made arbitrarily small by choosing $rn$ large enough, thanks to the definition of stabilization. Moreover, by the definition of $\ell_a$, see \eqref{elladef},
\[ I_{6rn}^{\mathrm{in}}(o) = \sum_{X_i \in X^\lambda \cap Q_{12rn\sqrt d}(o)} P_i \ell_{6rn}(|X_i|) \leq \ell(0) \sum_{X_i \in X^\lambda \cap Q_{12rn\sqrt d}(o)} P_i.  \]
In particular, using that the point process $\mathbf X^\lambda$ is independently marked with $P_i$ having marginal distribution $\mu$, and that $\Lambda$ is stationary with $\E[\Lambda(Q_1)]=1$, it follows that
\[ \mathbb E[I^{\mathrm{in}}_{6rn}(o)] \leq \ell(0) \lambda \E[P_o] \E[\Lambda(Q_{12rn\sqrt d})]=(12rn\sqrt d)^d \ell(0) \lambda \E[P_o]. \]
Thus, for any $n \geq 1$ and $r>0$, $\mathbb P(I_{6rn}^{\mathrm{in}}(o)>M)$ can be made arbitrarily small by choosing $M$ large enough, given that $\E[P_o]<\infty$. Thus, the statement of the lemma follows. 
%Proposition~\ref{prop-interferencecontrol} with $B_{r,n,M}(\cdot)$ replaced by $B_{r,n,M}^{\mathrm{in}}(\cdot)$ follows. 
\end{proof}

It remains to verify Proposition~\ref{prop-interferencecontrol} with $B_{r,n,M}(\cdot)$ replaced by $B_{r,n,M}^{\mathrm{out}}(\cdot)$ everywhere. More precisely, thanks to the exponential-moment and $b$-dependence assumption on $\Lambda$, the proof can be completed analogously to the proof of \cite[Proposition 3.3]{T18} starting from \cite[Equation (3.15)]{T18}, as soon as we have verified the following lemma. (\cite{T18} also assumed that $\ell(0) \leq 1$, but since $M$ can be made arbitrarily large in Proposition~\ref{prop-interferencecontrol}, $\ell(0) \leq 1$ can be assumed without loss of generality since for $\ell$ continuous, the function $\widetilde \ell=\ell/\ell(0)$ satisfies $\widetilde \ell(0)=1$, and for $a \geq 0$, we have $ \ell_{a}=\ell(0)\widetilde \ell_a$ and hence $I_a(x)=\ell(0)\sum_{i \in I} P_i \widetilde \ell_a(|x-X_i|)$.)
\begin{lem}\label{lemma-lookslikeT18}
Under the general assumptions of Theorem~\ref{thm-randomfadings} plus Part \eqref{second-fading}, there exists a constant $c_o=c_o(\mu,\ell)>0$ such that for all sufficiently small $s>0$, for all $\lambda>0$, $n \geq 1$ and $r>d_o$ with $rn>0$ sufficiently large and for all large enough $M>0$, for all $N\in\N$ and pairwise distinct $z_1,\ldots,z_N \in \Z^d$ we have
\begin{equation}\label{lookslikeT18}
\begin{split} 
\P(B^{\mathrm{out}}_{r,n,M}(z_1)^{\rm c}& \cap \ldots \cap B^{\mathrm{out}}_{r,n,M}(z_N)^{\rm c})\\ \leq & \E\Big[ \exp\Big( c_o \lambda s \sum_{i=1}^N \int_{\R^d\setminus Q_{12rn\sqrt d}(rnz_i)} \ell_{6rn}(|rnz_i-x|) \Lambda(\d x) \Big) \Big].
\end{split}
\end{equation}
\end{lem}
Indeed, the right-hand side of \eqref{lookslikeT18} is the same as the one of \cite[Equation (3.15)]{T18}, and the assumptions on $\Lambda$ in the two proofs are also the same.
%Now we conclude the proof of Proposition~\ref{prop-interferencecontrol} by verifying Lemma~\ref{lemma-lookslikeT18}.
\begin{proof}[Proof of Lemma~\ref{lemma-lookslikeT18}]
We start with an estimate originating from \cite[Section 3.2]{DF06}. By Markov's inequality, for any $s>0$,
\begin{align*}  
& \mathbb P( B^{\mathrm{out}}_{r,n,M}(z_1)^{\rm c} \cap \ldots \cap B^{\mathrm{out}}_{r,n,M}(z_N)^{\rm c})  = \mathbb P( I_{6rn}^{\mathrm{out}}(rnz_1)>M,\ldots,I_{6rn}^{\mathrm{out}}(rnz_N)>M ) 
\\  & \leq  \mathbb P \Big( \sum_{i=1}^N I_{6rn}^{\mathrm{out}}(rnz_i) > N M \Big)
\\ & \leq  \e^{-sNM} \E \Big[ \exp\Big( s \sum_{i=1}^N \sum_{X_k \in X^\lambda
\setminus Q_{12rn\sqrt d}(nz_i)
} P_k \ell_{6rn}(|rnz_i-X_k|) \Big) \Big].\numberthis\label{Csebisev} 
\end{align*}
The randomness of the power values $P_k$ prevents us from continuing the proof analogously to \cite{DF06,T18}. On the other hand, similarly to \cite[Section 4.3]{L19} in the Poisson case, we can argue as follows. According to the Marking Theorem~\cite[Section 5.2]{K93}, the independently marked CPP $\mathbf X^\lambda=\{ (X_i,P_i)\}_{i \in I}$ is a CPP in $\R^d \times [0,\infty)$ with directing measure $\Lambda \otimes \mu$, where we recall that $\mu=\P \circ P_o^{-1}$ is the distribution of $P_o$. Hence, applying the Laplace functional of a CPP (cf.~\cite[Sections~3.2,~6]{K93}) to the function $f \colon \R^d \times [0,\infty) \to [0,\infty)$, 
\[ f(x,p)=s\sum_{i=1}^N p \ell_{6rn}(|x-rnz_i|) \mathds 1 \{ x \in \R^d \setminus Q_{12rn\sqrt d}(rnz_i) \}, \] we obtain
\begin{align*}  
&\E \Big[ \exp\Big( s \sum_{i=1}^N \sum_{X_k \in X^\lambda
\setminus Q_{12rn\sqrt d}(rnz_i)
} P_k  \ell_{6rn}(|rnz_i-X_k|) \Big) \Big] \numberthis\label{CoxLaplace} \\ 
&= \E \Big[ \exp \Big( \lambda \int_{\R^d\setminus Q_{12rn\sqrt d}(rn z_i) }  \int_{0}^{\infty} \Big( \exp \Big( sp \sum_{i=1}^N \ell_{6rn}(|rnz_i-x|) 
\Big) -1 \Big) \mu(\d p)\Lambda(\d x)  \Big) \Big]. 
\end{align*}
Thanks to the exponential-moment assumption on $P_o$ from \eqref{second-fading}, the moment-generating function 
\[ \alpha \mapsto \E[\exp(\alpha P_o)]=\int_0^{\infty} \e^{\alpha p} \mu(\d p) \] is infinitely differentiable at $\alpha=0$ with first derivative $\int_0^{\infty} p \mu(\d p)=\E[P_o]<\infty$. 
Note that $\sum_{i=1}^N \ell_{6rn}(|rnz_i-x|) $ is uniformly bounded in $x\in \R^d$, $rn$, $N$ and pairwise distinct $z_1,\dots, z_N$, see~\cite[Lemma 3.6]{T18}. Consequently, for any $C>1$, the following holds for all sufficiently small $s>0$ (depending on $C$),
\[ \int_{0}^{\infty} \Big( \exp \Big( sp \sum_{i=1}^N\ell_{6rn}(|rnz_i-x|) \Big)-1 \Big) \mu(\d p) \leq Cs \E[P_o] \sum_{i=1}^N \ell_{6rn}(|rnz_i-x|). \numberthis\label{doubleexptrick2} \]
For such $s$, plugging \eqref{doubleexptrick2} back into \eqref{CoxLaplace}, starting from \eqref{Csebisev} we obtain
\begin{multline} \P(B^{\mathrm{out}}_{r,n,M}(z_1)^{\rm c} \cap \ldots \cap B^{\mathrm{out}}_{r,n,M}(z_N)^{\rm c}) \\ \leq \E\Big[ \exp\Big( C \E[P_o] \lambda s \sum_{i=1}^N \int_{\R^d\setminus Q_{12rn\sqrt d}(rnz_i)} \ell_{6rn}(|rnz_i-x|)\Lambda(\d x) \Big) \Big], \end{multline}
which is \eqref{lookslikeT18} with $c_o=C \E[P_o]$. With this we conclude the proof of the lemma.
\end{proof}

%========================================================

\subsubsection{Proof of Theorem~\ref{thm-randomfadings} Part~\eqref{first-fading}.}\label{sec-proof_thm1a}
Since $\Lambda$ is asymptotically essentially connected, the critical intensity $\lambda_{\rm c}(r)$ for the Poisson--Gilbert graph with connectivity threshold $r$ is finite, more precisely,  $\lambda_{\rm c}(r)<\infty$ holds for any $r>0$ according to \cite[Theorem 2.4]{HJC17}. Note further that the connection radii $(\rG^i)_{i \in I}$, defined in \eqref{r_B},  are bounded by $d_{\rm max}=\sup\{x\ge 0\colon x\in\supp(\ell)\}$. The proof of Theorem~\ref{thm-randomfadings} Part~\eqref{first-fading} can be obtained as an adaptation of the proof of Part~\eqref{second-fading} of the same theorem as follows.

First, one defines the subgraph of the SINR graph analogously to Step~\ref{step-subgraph} of the proof of Theorem~\ref{thm-randomfadings} Part~\eqref{second-fading}. Next, one takes the Step \ref{step-lattice} but for $r \in (d_o,d_{\rm max})$ arbitrary and fixed instead of letting $r \uparrow \infty$, and for $r_o(r)>r$ such that $r_o(r)$ still lies in the interval $(d_o,d_{\rm max})$ on which $\ell$ is strictly decreasing. This way, choosing $rn$ sufficiently large will be equivalent to choosing $n$ large enough (for fixed $r$). Further, one alters the choice of $\lambda(r)$: now, $\lambda(r)$ has to be chosen so large that the process of $(r,n)$-good sites percolates for some $n \geq 1$, which is possible for any fixed $r \in (d_o,d_{\rm max})$ since $\Lambda$ is asymptotically essentially connected, cf.~\cite[Section 5.2]{HJC17}. Next, Step~\ref{step-alsoSINRperc} is also applicable for all choices of the parameters where the underlying discrete model percolates. Finally, let us explain how to complete the proof of Proposition~\ref{prop-interferencecontrol} under the mere assumption that $\E[P_o]<\infty$. Since $\supp(\ell)$ is bounded, for all sufficiently large $n \geq 1$ the following holds for all $z \in \Z^d$
\[
\begin{aligned}
I_{6rn}(rnz) & = \sum_{X_i \in X^\lambda \cap Q_{6rn + 2d_{\rm max}}(rnz)} P_i \ell_{6rn}(|X_i-rnz|) \\ & \leq \sum_{X_i \in X^\lambda \cap Q_{12rn\sqrt d}(rnz)} P_i \ell_{6rn}(|X_i-rnz|)= I^{\mathrm{in}}_{6rn}(rnz). 
\end{aligned}\numberthis\label{allin}\]
Hence, it remains to control the inner part of the interference, which can be done analogously to Lemma~\ref{prop-interferencecontrolIN} once $\E[P_o]<\infty$, given that $\Lambda$ is stabilizing.  Hence, we conclude Theorem~\ref{thm-randomfadings} Part~\eqref{first-fading}.

\subsubsection{Proof of Theorem~\ref{thm-randomfadings} Part~\eqref{third-fading}.}\label{sec-proof_thm1c}
In the case when $\Lambda$ is only stabilizing and $P_{\sup}=\infty$, we observe that the proof of Theorem~\ref{thm-randomfadings} Part~\eqref{second-fading} stays valid if $\supp(\ell)$ is bounded but the following assumption holds: $\sup \supp(\ell)$ is sufficiently large such that
\[ \sup\supp(\ell) > \inf \{ r>0 \colon \exists n \geq 1 \text{ and }\lambda>0 \text{ such that $(r,n)$-good sites percolate} \}, \]
where the infimum is finite because $\Lambda$ is stabilizing. Indeed, in this situation, Lemma~\ref{lem-T18_1} (as in \cite[Section 3.2.2]{T18}) holds as well. 
Further, \eqref{allin} holds for all sufficiently large $n$ for all $z \in \Z^d$, and therefore one can complete the proof under the assumptions of Lemma~\ref{prop-interferencecontrolIN}, i.e., for $\Lambda$ stabilizing and $P_o$ such that $\E[P_o]<\infty$, without requiring $b$-dependence of $\Lambda$ for any $b>0$ or existence of exponential moments of $\Lambda(Q_1)$ or $P_o$.

\subsection{Proof of Proposition~\ref{prop-2degree1cluster}}\label{sec-twodegreeproof}
The strategy of the proof of the Proposition~\ref{prop-2degree1cluster} is the following. We first show that up to $\P$-null sets, clusters are either finite or infinite in both directions, i.e., they contain no vertex of degree 1 in case they are infinite, see Lemma~\ref{lemma-noheads} below. Next, we assume for a contradiction that there exists an infinite cluster with positive probability. Then, we introduce a procedure that removes points from the infinite cluster that is closest to the origin in a certain sense. Thanks to elementary properties of the SINR graph, in the resulting configuration, 
the infinite cluster still remains infinite, but it contains a vertex of degree 1. Hence, the probability that the process takes values in the set of the resulting configurations is zero. What remains to show afterwards is that also the probability that the process takes place in the set of original configurations is zero, which leads to the desired contradiction. At this point it will be useful to compare the resulting configuration with an independent thinning of the original configuration in a certain ball, and this is where we make use of the fact that the underlying point process is a stationary CPP. 

%Throughout the proof, we will assume that $\P(P_o>0)=1$; indeed, vertices $X_i$ of the Cox--SINR graph with $P_i=0$ are disconnected in the graph, and thus if $g_{(\g)}(\mathbf X^\lambda)$ percolates, then so does $g_{(\g)}(\mathbf X^{\lambda,+})$, where $X^{\lambda,+}=\{ X_i \in X^\lambda \colon P_i>0 \}$ is an independent thinning of $X^\lambda$ with survival parameter $p=\P(P_i>0)$.
We assume throughout the proof that $\gamma \geq 1/(2\tau)$, so that degrees in $\Gg$ are bounded by two, and that $X^\lambda$ is nonequidistant (for all $\lambda>0$). We can also assume that $\P(P_o>0)>0$ in what follows, since otherwise the statement is trivially true. 
We start the proof with the following lemma, which excludes infinite paths that have an endpoint in case the degrees are bounded by two, in a substantially more general setting.
\begin{lem}\label{lemma-noheads}
Let $g(\X)$ be a random graph based on a stationary marked point process $\X=\{(X_i,P_i)\}_{i\in I}$ with values in $\R^d \times Z$, where the mark space $(Z,\mathcal Z)$ is an arbitrary measurable space, $X=\{X_i\}_{i\in I}$ is the vertex set, and such that the degree of all $X_i\in X$, $\deg(X_i)$, is bounded by $2$, almost surely. Let $X$ have a finite intensity and consider the point process of degree-one points in infinite clusters
\[ \mathcal X_0=\sum_{i\in I}\delta_{X_i}\one\{\deg(X_i)=1,\, X_i \text{ is part of an infinite cluster in }g(\X)\}.\]
Then, $\P(\mathcal  X_0(\R^d)=0)=1$. 
%all points contained in an infinite cluster of $\Gg$ have degree equal to two, almost surely.
\end{lem}
We will apply this lemma to the SINR graph $g(\X)=\Gg$ with $\lambda$ arbitrary and $\gamma \geq 1/(2\tau)$ and $Z=[0,\infty)$. The proof is based on a variant of the \emph{mass-transport principle} (cf.~\cite[Section 4.2]{CDS20} for instance). 
%To be more precise, the proof goes as follows.
\begin{proof}[Proof of Lemma~\ref{lemma-noheads}.]
First, using the union bound and stationarity, it is enough to show that $\E[\mathcal  X_0(Q_1)]=0$. Let us define the point process of points in infinite clusters that are at distance equal to $k\in\N_o$ from a point in $\mathcal X_0$, 
$$
\mathcal X_k=\sum_{i\in I}\delta_{X_i}\one\{X_i\text{ is part of an infinite cluster and has graph distance }k\text{ from }\mathcal X_0\}.
$$
Thanks to the degree bound, every infinite cluster has at most one point in $\mathcal X_0$ and $\E[\mathcal X_k(Q_1)]=\E[\mathcal X_0(Q_1)]$, for all $k\in\N_o$, by stationarity. However, $\sum_{k\ge 0}\E[\mathcal X_k(Q_1)]\le \E[X(Q_1)]<\infty$ and thus $\E[\mathcal X_0(Q_1)]=0$. 
\end{proof}

Let us denote by $(\mathcal C_i)_{0\le i< L}$ the $L$-many infinite clusters in $\Gg$, where $L\in \N\cup\{ \infty\}$. For the proof of Proposition~\ref{prop-2degree1cluster}, it then suffices to show that
\[ \P(L\ge 1)=0. \numberthis\label{onotinaninfinitecluster} \] 

%Let us write $P=\P\circ \X^\lambda^{-1}$ for the distribution of $\X^\lambda$. 
%
We view $\X^\lambda$ as the canonical process $\X^\lambda(\bs\omega)=\bs\omega$ on the set $\mathbf N$ of marked point configurations $\bs \omega$ in $\R^d  \times [0,\infty)$ such that $\omega=\{ x_i \colon (x_i,p_i) \in \bs \omega \}$ is an infinite, locally-finite, nonequidistant point configuration on $\R^d$.  The set of such point configurations $\omega$ will be denoted by $\rm N$. Note that $\mathbf N$ and ${\rm N}$ are equipped with the corresponding  evaluation $\sigma$-fields. 
%The main advantage of using the configuration space $\mathbf N$ that we can completely disregard certain $\P$-null sets such as the ones that are not nonequidistant.  %Note that if $\bs \omega,\bs \omega' \in \mathbf N$ are such that $\bs \omega \subseteq \bs \omega'$, then for any $x,y \in \omega$ such that $\SINR((x,p),(y,q),\bs\omega')>\tau$, we also have $\SINR((x,p),(y,q),\bs\omega)>\tau$. Hence, $g_{k,f}(\bs\omega)$ contains all edges of $g_{k,f}(\bs\omega')$ that connect two points of $\omega$.

Now we introduce an ordering in $\R^d\times [0,\infty)$, which orders the points of the set according to the received signal power at a given point $y \in \R^d$ (or equivalently, according to the received SINR values $\SINR(\cdot,y,\bs \omega)$ at $y$). 
\begin{defn}\label{defn-closer}
Let $(x,p),(z,r)\in\R^d \times [0,\infty)$ and $y \in \R^d$. We say that \emph{$(x,p)$ transmits a stronger signal to $y$ than $(z,r)$ does} if one of the following conditions is satisfied:
\begin{enumerate}
    \item\label{signalordering} $p \ell(| x-y |) > r \ell(| z-y |)$, or
    \item\label{tiebreaking} $p \ell(| x-y |) = r \ell(| z-y |)$ and $| x-y |<| z-y|$.
\end{enumerate}
\end{defn}
When talking about the marked CPP $\X^\lambda$, we will always assume that transmitters are associated with their own transmitted signal powers, and hence we will say ``$X_i$ transmits a stronger signal to $X_j$ than $X_l$ does'' instead of ``$(X_i,P_i)$ transmits a stronger signal to $X_j$ than $(X_l,P_l)$ does'', for any $i,j,l \in I$ such that $i \neq j$ and $l \neq j$. Then, it is easy to see that for $\X^\lambda$ such that $X^\lambda$ is nonequidistant, almost surely the following holds. For all $i \in I$, the relation ``$X_i$ transmits a stronger signal to $X_j$ than $X_l$ does'' is a total ordering (i.e., irreflexive, antisymmetric and transitive, with any two elements being comparable) on the set $\{ (i,l) \in I^2 \colon i \neq j \text{ and } l \neq j \}$, which we call the \emph{ordering of signal-weighted distance} from receiver $X_i$. This fact indeed relies on the tiebreaking mechanism \eqref{tiebreaking}: e.g., if $\ell$ is constant on some interval (which is possible under the assumption of Proposition~\ref{prop-2degree1cluster} and even under the stronger assumption of Theorem~\ref{thm-randomfadings}), then \eqref{signalordering} does not define a total ordering on its own.

For $\bs\omega \in \mathbf N$ and $x_o \in \omega$, we can consider the vector $\mathbf V(x_o,\bs \omega)=(\mathbf \NNN_n(x_o,\bs \omega))_{n \in \N_0}$ of the marked points of $\bs \omega$ ordered increasingly according to signal-weighted distance from receiver $x_o$.
Then, we define $\NN_i(x_o,\bs \omega)$ as the first component of $\NNN_i(x_o,\bs \omega)$, which we call the \emph{$i$-th nearest neighbor of $x_o$ in signal-weighted order}. In particular, $\NN_0(x_o,\bs \omega)=x_o$. Note that if the distribution $\mu$ is concentrated in one point $p>0$, i.e., $\mu=\delta_p$, then the $i$-th nearest neighbor of $x_o$ in signal-weighted order is just the $i$-th nearest neighbor of $x_o$ with respect to Euclidean distance.  
%We will use the notation $\NN_i(\bs\omega)=\NN_i(o,\bs \omega)$ and also write $\NNN_i(x_
%We will use the notation $\NN_i(\bs\omega)=\NN_i(o,\bs \omega)$ and also write $\NNN_i(x_o,\bs \omega)$ for the $i$-th entry of $\NNN(x_o,\bs \omega)$.

Now, if $x_o$ has degree two in $G_\g(\X^\la(\bs\omega))$, then $x_o$ must be connected by an edge to both $\NN_1(x_o,\bs\omega)$ and $\NN_2(x_o,\bs\omega)$ since the degree bound applies already for the edges towards $x_o$. Moreover, both $\NN_1(x_o,\bs\omega)$ and $\NN_2(x_o,\bs\omega)$ must also have $x_o$ as one of their first two nearest neighbors in signal-weighted order, that is, 
\[ x_o \in\big \{ \NN_1(\NN_i(x_o,\bs\omega),\bs\omega), \NN_2(\NN_i(x_o,\bs\omega),\bs\omega) \big\}, \]
for all $i\in\{1,2\}$. 
These signal-weighted nearest neighbor relations hold almost surely, in particular for every nonequidistant configuration $\bs\omega$. The goal of using the configuration space $\mathbf N$ is to entirely exclude configurations that offend the degree bound or the signal-weighted nearest neighbor relations or are not nonequidistant.

In the event $\{L\ge 1\}$, let $\mathbf Z=(Z,R)$ denote the closest point to the origin that has degree two and is contained in an infinite cluster. Without loss of generality, we will assume that this cluster is always equal to $\mathcal C_0$. Now, Proposition~\ref{prop-2degree1cluster} immediately follows once we have verified the following proposition. 
\begin{prop}\label{lemma-inprop2deg}
Consider the event $\{L\ge 1\}$ and define the random variable
\[ I = \inf \{ i \geq 3 \colon \NN_i(Z,\X^\la) \in \mathcal C_0 \}. \]
Then, under the assumptions of Proposition~\ref{prop-2degree1cluster}, for any $i \geq 3$, we have 
\[ \P(\{ L\ge 1 \}\cap \{ I=i \})=0. \numberthis\label{there'snobody} \]
\end{prop}

\begin{proof}[Proof of Proposition~\ref{prop-2degree1cluster}]
Using a union bound and noting that $\{ L\ge 1  \}\subset \{I<\infty\}$, Proposition~\ref{lemma-inprop2deg} implies $\P(L\ge 1)=0$, which is \eqref{onotinaninfinitecluster}, and thus the proof of Proposition~\ref{prop-2degree1cluster} is finished.
%In order to conclude that $\P(A)=0$, it remains to verify the following lemma.
%\begin{lem}\label{lemma-omegaomegai}
%Having proven Proposition~\ref{lemma-inprop2deg}, Proposition~\ref{prop-2degree1cluster} also follows.
%\end{lem}
\end{proof}

\begin{proof}[Proof of Proposition~\ref{lemma-inprop2deg}]
%Let us drop the index $j$ in $\mathcal C_j$ and simply write $\mathcal C$, and the same for $\mathbf Z_j$ and $I_j$. 
For $\bs\omega\in \{L\ge 1\}$, by definition, we have that $Z(\bs\omega)$ is connected by an edge both to $\NN_1(Z(\bs\omega), \bs\omega)$ and $\NN_2(Z(\bs\omega), \bs\omega)$ in $G_\g(\X^\la(\bs\omega))$. Further, thanks to the degree bound of two, in the event $\{L\ge 1\}$, $\NN_1(Z(\bs\omega), \bs\omega)$ and $\NN_2(Z(\bs\omega), \bs\omega)$ have no further joint neighbor in $G_\g(\X^\la(\bs\omega))$ since otherwise $\mathcal C_0(\bs\omega)$ has a loop and cannot be infinite by the degree bound. 
%the cluster containing them would contain a loop of order 4 and thus could not be infinite. 
This way, for any $i \geq 3$, there exists $l \in \{1,2\}$ such that $\NN_i(Z(\bs\omega), \bs\omega)$ and $\NN_l(Z(\bs\omega), \bs\omega)$ are not connected by an edge in $G_\g(\X^\la(\bs\omega))$. Let us denote the corresponding $\NN_l(Z(\bs\omega), \bs\omega)$ by $M_i(\bs\omega)$, and define $M_i(\bs\omega)=\NN_1(Z(\bs\omega), \bs\omega)$ if neither $\NN_1(Z(\bs\omega), \bs\omega)$ nor $\NN_2(Z(\bs\omega), \bs\omega)$ is connected to $\NN_i(Z(\bs\omega), \bs\omega)$ by an edge. The element of $\{ \NN_1(Z(\bs\omega), \bs\omega),\NN_2( Z(\bs\omega), \bs\omega) \}$ not being equal to $M_i(\bs\omega)$ is denoted by $N_i(\bs\omega)$.
We will write $Q$ for the signal power transmitted by $M_i(\bs\omega)$. 

Let us fix $i \geq 3$. Let $\bs \omega \in \{L\ge 1\}$ be such that $I(\bs \omega)=i$. Let us define a thinned configuration
\begin{align*} 
\bs \omega^i =\bs \omega\setminus \{(M_i(\bs\omega),Q), \NNN_3(Z(\bs\omega),\bs\omega), \dots, \NNN_{i-1}(Z(\bs\omega),\bs \omega)\}.
\end{align*}
%We abbreviate the thinning by $h_i(\bs\omega)=\bs\omega^i$. 
%We claim that $\bs\omega^i \in \mathbf N$. Indeed, for any $\omega \in {\rm N}$, $\omega$ minus a finite set of points in $\R^d$ is an element of ${\rm N}$, which holds in particular for $\omega^i$. %Further, the power values of $\bs\omega^i$ are also chosen from the suitable space, hence the claim follows.
%Note that the thinned configurations is designed in such a way that 

We claim for $\P$-almost all $\bs \omega \in \{L\ge 1\}\cap \{I=i\}$ also $\bs \omega^i \in \{L\ge 1\}$. For this, first note that the removal of finitely many points and their associated edges from an infinite cluster does not change the property of the cluster to be infinite. However, the removal of points can still change the edge structure of the remaining points. In order to exclude this, we can use the following fundamental property of the SINR graph. Assume that $\bs \omega,\bs \omega'$ are elements of $\mathbf N$ such that $\bs \omega \subseteq \bs \omega'$. Then, for all $x,y \in \omega$, if $\SINR(x,y,\bs \omega') >\tau$, then $\SINR(x,y,\bs \omega) >\tau$, which is clear from \eqref{SINRproof}. In words, if we remove some vertices from an SINR graph, then edges of the SINR graphs between the remaining points stay preserved.

Then, our next claim is that for $\bs\omega\in \{L\ge 1\}\cap \{I=i\}$, we have that $\bs\omega^i$ is contained in
\[ B = \{\bs\eta\colon  L(\bs\eta)\ge 1 \text{ and }\mathcal C_0(\bs\eta) \text{ contains a point of degree one}\} \subset \{L\ge 1\}.  \]
The proof of this claim in the simplest case $i=3$ is illustrated in Figure~\ref{fig-NN3}.
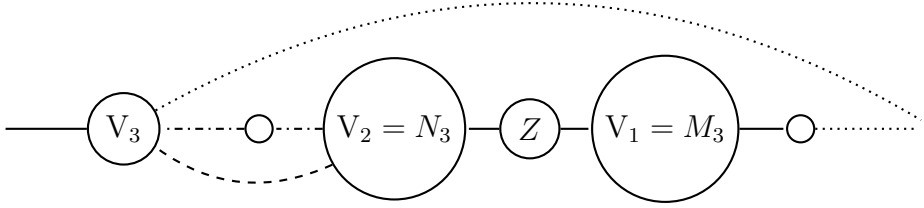
\begin{figure}[!htpb]
\centering
\begin{tikzpicture}[scale=3,
        shorten >=1pt, auto, thick,
        node distance=1.8cm,
    main node/.style={circle,draw,font=\sffamily\large\bfseries}
                            ]

			\node[main node] (1) {$Z$};
			\node[main node] (2) [right of=1] {$\NN_1=M_3$};
			\node[main node] (3) [left of=1] {$\NN_2=N_3$};
			\node[main node] (5) [right of=2] {};
			\node[main node,color=white] (8) [right of=5] {};
			\node[main node] (6) [left of=3] {};
			\node[main node] (4) [left of=6] {$\NN_3$};
			\node[main node,color=white] (7) [left of=4] {};

			\path[every node/.style={font=\sffamily\small}]
			  (1) edge (2)
			      edge (3)   
			  (6) edge [dashdotted] (4)
			  (4) edge (7)
			  (2) edge (5)
			  (5) edge [dotted] (8)
			  (3) edge [dashdotted] (6)
			  (3) edge [dashed,bend left] (4)
			  (4) edge [dotted, bend left] (8);
\end{tikzpicture}
\caption{A visualization of the case $I(\bs\omega)=3$ for some realization $\bs\omega\in \{L \geq 1\}$. $\NN_3=\NN_3(Z(\bs \omega),\bs\omega)$ is contained in the infinite cluster $\mathcal C_0=\mathcal C_0(\bs\omega)$ of the SINR graph $G_{\g}(\bs \omega)$ including $Z=Z(\bs \omega)$, and it is not a neighbor of $M_3=M_3(\bs\omega)$, which in this example equals $\NN_1=\NN_1(Z(\bs\omega),\bs\omega)$, whereas $\NN_2=\NN_2(Z(\bs \omega),\bs\omega)=N_3=N_3(\bs\omega)$. Hence, if $\NN_3$ has degree two in $\mathcal C_0$,  then there are various possibilities respecting the degree bound of two to connect $\NN_3$ to $\mathcal C_0$ so that it is not connected to $M_3$ by an edge. $\NN_3$ can either be a direct neighbor of $\NN_2$ (see dashed line) or a later point of the path from $Z$ to infinity starting with the edge from $Z$ to $\NN_2$ (dash-dotted lines) or a non-direct neighbor of $\NN_1$ on the path from $Z$ to infinity starting with the edge from $Z$ to $\NN_1$ (dotted lines).
Now, removing $M_3$ from the realization, both edges adjacent to $\NN_3$ are preserved. Also all edges from $Z$ to infinity starting with the edge from $Z$ to $\NN_2$ are preserved, hence $Z$ is still contained in an infinite cluster, but the edge from $Z$ to $\NN_1$ is removed. In the resulting new configuration, the second-nearest neighbor of $Z$ in signal-weighted order is $\NN_3$, and hence this is the only point of the configuration that could be connected to $Z$ by an edge. But $\NN_3$ still cannot have degree 3 or more, hence it cannot be connected to $Z$, which implies that in the new configuration $Z$ is in an infinite cluster containing a point of degree one.}
\label{fig-NN3}
\end{figure}
%In order to see this, note that for $\P$-almost all $\bs \omega \in \{L\ge 1\}\cap \{I=i\}$, the following two conditions are both satisfied.
%\begin{enumerate}[(i)]
%\item\label{first-goodconf} There are precisely two edge-disjoint infinite paths in $\Gg(\bs\omega)$ starting from $Z(\bs\omega)$. Hence, in particular, at least one of these paths does not pass through $M_i(\bs\omega)$ and
%\item\label{second-goodconf} $\NN_3(Z(\bs\omega),\bs\omega),\ldots,\NN_{i-1}(Z(\bs\omega),\bs\omega) \notin \mathcal C_0(\bs\omega)$ by the definition of $I$ and the fact that $I(\bs \omega)=i$. 
%\end{enumerate}
For general $i\ge 3$, recall that $Z$ cannot have degree higher than two in $G_\g(\X^\la(\bs\omega^i))$, whereas it has degree at least one and its cluster $\mathcal C_0(\bs\omega^i)$ is infinite in $G_\g(\X^\la(\bs\omega^i))$. Note also that the edge between $Z(\bs\omega)$ and $N_i(\bs\omega)$ still exists in $G_\g(\X^\la(\bs\omega^i))$. Further, if $Z(\bs\omega)$ has degree two in $G_\g(\X^\la(\bs\omega^i))$, then it is connected to the second-nearest neighbor of $Z(\bs\omega)$ in signal-weighted order in $\bs\omega^i$, which is $\NN_2(Z(\bs\omega),\bs\omega^i)=\NN_i(Z(\bs\omega),\bs\omega)$, whereas $\NN_1(Z(\bs\omega),\bs\omega^i)=N_i(\bs\omega)$. Now, since $\bs\omega \notin B$, $\bs\omega \in \{L\ge 1\}$ and $\NN_i(Z(\bs\omega),\bs\omega) \in \mathcal C_0(\bs\omega)$, it follows that $\NN_i(Z(\bs\omega),\bs\omega)$ has degree equal to two in $G_\g(\X^\la(\bs\omega))$. Further, it is neither connected to $M_i(\bs\omega)$ by an edge nor to $Z(\bs\omega)$ in this graph. Hence, both edges adjacent to $\NN_i(Z(\bs\omega),\bs\omega)$ also exist in $G_\g(\X^\la(\bs\omega^i))$. But since $\NN_i(Z(\bs\omega),\bs\omega)$ has degree at most two in $G_\g(\X^\la(\bs\omega^i))$, it follows that $Z(\bs\omega)$ and $\NN_i(Z(\bs\omega),\bs\omega)$ are not connected by an edge in this graph. Hence, $\bs\omega^i \in B$, which implies the claim.

Note that by Lemma~\ref{lemma-noheads}, the set $B$ is a $\P$-null set, i.e.,
%In other words, since $\mathbf X^\lambda$ is the canonical process on $\mathbf N$ with distribution $\P$, 
\[ \P(\{ \bs \omega^i \colon \bs \omega \in \{L\ge1\}\cap \{ I=i \} \})=0.\numberthis\label{isunlikely} \] 
This implies \eqref{there'snobody} and concludes the proof of Proposition~\ref{lemma-inprop2deg} as soon as the following lemma is verified. 
\begin{lem}\label{lemma-omegaomegai} 
Under the assumptions of Proposition~\ref{prop-2degree1cluster}, for any $i \geq 3$, $\P(\{L\ge1\}\cap \{ I=i \})>0$ implies $\P( \{ \bs \omega^i \colon \bs \omega \in \{L\ge1\}\cap \{ I=i \} \}) >0$.
\end{lem}
By Lemma~\ref{lemma-omegaomegai}, where we show that if the collection of thinned configurations is contained in a $\P$-null set, also the non-thinned configurations form a $\P$-null set, we see that \eqref{isunlikely} implies \eqref{there'snobody}, which concludes the proof of Proposition~\ref{lemma-inprop2deg}.
\end{proof}

\begin{proof}[Proof of Lemma~\ref{lemma-omegaomegai}]
Let us fix $i \geq 3$ and assume that $\P(\{L\ge1\}\cap \{ I=i \})>0$. Then, by continuity of measures, there exists $K>0$ such that 
\[ \P\big( \big\{ \bs \omega \in \{L\ge1\}\cap\{I=i\} \colon \NN_j(Z(\bs\omega), \bs \omega) \in B_K(o),~\forall j \in \{ 1,\ldots,i\}\big \}\big)>0, \]
where $B_K(o)$ denotes the open Euclidean ball of radius $K$ in $\R^d$. Hence, there exists $n \geq i$ such that $\P(C_{i,K,n})>0$, where 
\[
\begin{aligned}
 C_{i,K,n} =\big\{ \bs \omega \in \{L\ge1\}\cap \{I=i\} \colon &\#\big( \omega \cap B_K(o) \big)=n+1 \\
 &\text{and }\NN_j(Z(\bs\omega),\bs \omega) \in B_K(o), \forall j \in \{ 1,\ldots,i\}\big\}. \end{aligned}  \]
Conditional on the event $C_{i,K,n}$, the marked point process $(\mathbf X^\lambda \setminus \{\mathbf Z \})\cap B_K(o)$ has precisely $n$ points $\mathbf X_1,\dots \mathbf X_n$.

Now, for some fixed $q \in (0,1)$, we can represent $\mathbf X^\lambda$ as $\X^{\lambda,1} \cup \X^{\lambda,2}$ as follows.  For $K>0$, let $B_K(o)$ denote the open $\ell^2$-ball of radius $K$ around $o$. Let $\X^{\lambda,1}$ be given as the  union of $\mathbf X^\lambda \setminus (B_K(o)\times [0,\infty))$ and the independent thinning of $\mathbf X^\lambda \cap (B_K(o)\times [0,\infty))$ with survival probability $q$, and let $\X^{\lambda,2}$ be the complementary thinning. That is, conditional on $\mathbf X^\lambda$, $\X^{\lambda,1} \cap (B_K(o)\times [0,\infty))$ contains each point of $\mathbf X^\lambda \cap (B_K(o)\times [0,\infty))$ with probability $q$ independent of the other points of this point process, and it contains no other points. Note further that $\X^{\lambda,2}$ and  $\X^{\lambda,1} \cap (B_K(o)\times [0,\infty))$ are independent thinnings of $\mathbf X^\lambda \cap (B_K(o)\times [0,\infty))$ with survival probability $1-q$ respectively $q$, moreover, $\X^{\lambda,1}=\mathbf X^\lambda \setminus \X^{\lambda,2}$, $\X^{\lambda,2} \setminus (B_K(o)\times [0,\infty))=\emptyset$, and $\X^{\lambda,1} \setminus (B_K(o)\times [0,\infty))=\mathbf X^\lambda \setminus (B_K(o)\times [0,\infty))$. 
In order to provide a precise construction of the thinned processes, we choose a sequence $(J_m)_{m \in\N}$ of i.i.d.~Bernoulli random variables with parameter $q$ that is independent of $\mathbf X^\lambda$, and given the realization $\bs \omega=\mathbf X^\lambda(\bs \omega) = (\NNN_i(Z(\bs\omega),\bs \omega))_{i \in \N_0}$, the realizations of $\X^{\lambda,1}(\bs \omega)$ and $\X^{\lambda,2}(\bs \omega)$ are defined as follows, depending also on $(J_m)_{m \in\N}$:
\[\begin{aligned}\X^{\lambda,1}(\bs \omega) & =\X^{\lambda,1}(\bs \omega,(J_m)_{m \in \N}) = \{ \NNN_m(Z(\bs\omega),\bs \omega) \colon J_m = 1, \NNN_m(Z(\bs\omega),\bs \omega) \in B_K(o) \}\\ 
& \quad \cup \{ \mathbf Z(\bs\omega)  \} \cup  \{ \mathbf \NNN_m(Z(\bs\omega),\bs \omega) \colon \NNN_m(Z(\bs\omega), \bs \omega) \in \R^d \setminus B_K(o) \} \end{aligned} \]
and
\[\X^{\lambda,2}(\bs \omega)=\X^{\lambda,2}(\bs \omega,(J_m)_{m \in \N}) = \{ \NNN_m(Z(\bs\omega), \bs \omega) \colon J_m = 0, \NNN_m(Z(\bs\omega),\bs \omega) \in B_K(o) \}. \]
It is clear that the projections $X^{\lambda,1},X^{\lambda,2}$ of $\X^{\lambda,1}$ respectively $\X^{\lambda,2}$ to the $\R^d$-coordinate are nonequidistant, further, $\X^{\lambda,1}$ can be represented as a random variable with values in $\mathbf N$, defined on an enlarged probability space  $(\Omega',\mathcal F',\mathbb P')$ governing both the point process $\X^\lambda$ and the sequence $(J_m)_{m \in \N}$. In particular, $\P'(\X^\lambda \in \cdot) = \P(\X^\lambda \in \cdot)$.

The next property of stationary and nonequidistant CPPs is crucial for completing the proof of Lemma~\ref{lemma-omegaomegai}.
\begin{lem}\label{claim-stablethinning}
Let $\Lambda$ be stationary and nonequidistant. Then, for any $K>0$, the law of $\X^{\lambda,1}$ is absolutely continuous with respect to the one of $\X^\lambda$ for any $K>0$.
\end{lem}
To be more precise, the absolute continuity is meant in this lemma in the following way, with respect to the probability space $(\Omega',\mathcal F',\mathbb P')$ on which $\X^{\lambda,1}$ and $\X^\lambda$ are jointly defined with $\mathbb P'(\X^\lambda \in \cdot)=\mathbb P(\X^\lambda \in \cdot)$. Let $G \in \mathcal F'$ be any event such that $\P'(\X^{\lambda,1} \in G)>0$, then we have $\P'(\X^\lambda \in G)>0$. 
\begin{proof}[Proof of Lemma~\ref{claim-stablethinning}.]
 Let $F$ be an element of the evaluation $\sigma$-algebra of $\mathbf N$ such that $\P'(\X^{\lambda,1} \in F)>0$. We have to show that then also $\P(\mathbf X^\lambda \in F)>0$. Under the assumption that $\P'(\X^{\lambda,1} \in F)>0$, by continuity of measures, we can find $K,l \in \N$ such that
\[ \eps:=\P'(\X^{\lambda,1} \in F, \# (\X^{\lambda,1} \cap (B_K(o) \times [0,\infty))) = l) >0. \numberthis\label{eps1prob}\]
In other words, we have $0<\eps = \P'(\X^{\lambda,1} \in G)$ where $G = \{ \bs \omega \in F \colon \# (\bs \omega \cap ( B_K(o) \times [0,\infty)))=l \}$.
Thus,
\[
\begin{aligned}
 \P(\mathbf X^\lambda \in F)  &\geq \P'(\mathbf X^\lambda \in G, \X^{\lambda,1} = \mathbf X) 
 \geq \P'(\X^{\lambda,1} \in G)\P'(\X^{\lambda,1} = \mathbf X^\lambda | \X^{\lambda,1} \in G) \\ &= \eps \P'(\X^{\lambda,1} = \mathbf X^\lambda | \X^{\lambda,1} \in G),
\end{aligned}
\] 
 and further,
\[ 
\begin{aligned}
 \P'(\X^{\lambda,1} = \mathbf X^\lambda | \X^{\lambda,1} \in G) & %=  \frac{\P'(\X^{\lambda,1} = \mathbf X, \X^{\lambda,1} \in G)}{\P'( \X^{\lambda,1} \in G)}  
 \geq \P'(\X^{\lambda,1} = \mathbf X, \X^{\lambda,1} \in G)
 = \P'(\X^{\lambda,2} = \emptyset, \X^{\lambda,1} \in G).
\end{aligned} \numberthis\label{probofcond} \]
According to~\eqref{eps1prob}, we have
\[ 0<\eps = \P'(\X^{\lambda,1} \in G) = \sum_{n=0}^{\infty}a_n,\]
where $a_n=\E'\big[ \P'(\X^{\lambda,1} \in G | \Lambda) \mathds 1 \{ \Lambda(B_K(o)) \in [n,n+1) \} \big]$, and thus there exists $m\in\N_0$ with $a_m>0$. 
Now, conditional on $\Lambda$, $\X^{\lambda,1}$ is an i.i.d.~marked PPP, and hence a PPP on $\R^d \times [0,\infty)$, which also implies that the complementary thinnings $\X^{\lambda,1}$ and $\X^{\lambda,2}$ are independent given $\Lambda$, see \cite[Colouring Theorem and Marking Theorem]{K93}. Hence, we obtain
\[
\begin{aligned} 
\P'(\X^{\lambda,2} = \emptyset, &\X^{\lambda,1} \in G) = \E'\big[\P'(\X^{\lambda,2} = \emptyset|\Lambda) \P'(\X^{\lambda,1} \in G|\Lambda) \big] \\
%& = \sum_{n=0}^{\infty}\E'\big[ \e^{-\lambda(1-p)\Lambda(B_K(o))} \P'(\X^{\lambda,1} \in G, \Lambda(B_K(o)) \in [n,n+1) |\Lambda) \big]
& = \sum_{n=0}^{\infty} \E'\big[   \e^{-(1-q)\Lambda(B_K(o))}\P'(\X^{\lambda,1} \in G |\Lambda)  \mathds 1\{ \Lambda(B_K(o)) \in [n,n+1)  \} \big]
\\ & \geq \sum_{n=0}^{\infty} \e^{-(1-q)(n+1)} a_n\ge \e^{-(1-q)(m+1)} a_m>0,
\end{aligned}
\]
which verifies the lemma that the distribution of $\X^{\lambda,1}$ is absolutely continuous with respect to the one of $\X^\lambda$.
\end{proof}
Given Lemma~\ref{claim-stablethinning}, we now finish the proof of Lemma~\ref{lemma-omegaomegai}.
%\begin{comment}
%{We now claim that under $\mathbb P'$, the distribution of $\X^\lambda^{*}$ is absolutely continuous with respect to the one of $\X^{\lambda,1}$. Indeed, let $G$ be any element of the evaluation $\sigma$-algebra of $\mathbf N$ and assume that $\P(\X^\lambda \in G)>0$. In order to verify the claim, our goal is now to show that $\P'(\X^{\lambda,1} \in G)>0$. Then, since $\X^\lambda$ is locally finite and contains $(o,p_o)$ almost surely, there exists $m \in \N$ such that 
%\[ \eps_0:=\P(\X^\lambda \in G, \# (\X^\lambda \cap (B_K(o) \times [0,\infty))) = m ) >0. \]
%Thus, we can estimate
%\[ 
%\begin{aligned}
%\P'&(\X^{\lambda,1} \in G) \geq \P(\X^{\lambda,1} \in G, \# (\X^{\lambda,1} \cap (B_K(o) \times [0,\infty)) ) = m ) \\ &\geq \P'(\X^{\lambda,1} \in G, \# (\X^{\lambda,1} \cap (B_K(o) \times [0,\infty)) ) = m ) \\ &  \geq \P'(\X^\lambda^{*} \in G, \# (\X^\lambda^{*} \cap (B_K(o) \times [0,\infty)) ) = m, \X^{\lambda,1} = \X^\lambda) \\ & \geq \P(\X^\lambda^{*} \in G, \# (\X^\lambda^{*} \cap (B_K(o) \times [0,\infty)) ) = m) \P'( \X^{\lambda,1} = \X^\lambda | \X^\lambda^{*} \in G, \# (\X^\lambda^{*} \cap (B_K(o) \times [0,\infty)) ) = m) \\ & = \eps_0 p^{m-1} >0,
%\end{aligned} \numberthis\label{firstabscont}
%\]
%which implies to claim.}
%\end{comment}
Thanks to the assumption that $\P(C_{i,K,n})>0$ and using the definition of $\X^{\lambda,1}$,
\[
\begin{aligned}
 \mathbb P' \big( \X^{\lambda,1} \in& \{ \bs \omega^i \colon \bs \omega \in \{L\ge 1\} \cap \{ I=i \} \}\big)  \geq \mathbb P' \big( \X^{\lambda,1} \in \{ \bs \omega^i \colon \bs \omega \in C_{i,K,n} \} \big) \\ & \geq \mathbb P' \big( \X^{\lambda,1} \in \{ \bs \omega^i \colon\bs \omega \in C_{i,K,n}\}, \X^\lambda \in C_{i,K,n}  \big)
 \\ & = \P(C_{i,K,n}) \mathbb P'\big(\X^{\lambda,1} \in \{\bs \omega^i \colon \bs \omega \in C_{i,K,n} \}| \X^\lambda \in C_{i,K,n} \big) \\
 & = \P(C_{i,K,n}) q^{n-i+2} (1-q)^{i-2}>0. 
\end{aligned} \numberthis\label{thinnedestimate}
\]
Finally, by Lemma~\ref{claim-stablethinning}, under $\mathbb P'$ the distribution of $\X^{\lambda,1}$ is absolutely continuous with respect to the one of $\X^\lambda$. Hence, it follows from~\eqref{thinnedestimate} that
\[  \mathbb P \big( \X^\lambda \in \{ \bs \omega^i \colon \bs \omega \in \{L\ge 1\} \cap \{ I=i \} \} \big)>0, \]
which implies the lemma.
\end{proof}

\subsection{Proof of Theorem~\ref{thm-PPPlambdaequality}}\label{sec-lambda=lambdaproof}
This proof is similar to the one of Theorem~\ref{thm-randomfadings} Part~\eqref{second-fading} but simpler. The new proof ingredient that we use here is the strong connectivity of \emph{any} supercritical Poisson--Boolean model \cite[Theorems 2 and 5]{PP96} in case $d\geq 2$, which allows us to improve the result that $\lacc<\infty$ to $\lacc=\lambda_{\rm c}(r_{\rm B})$. First we introduce an adequate discrete percolation model and then we control the interferences. 

%We use the notation of Section~\ref{sec-IntroMainRes} for the setting of Theorem~\ref{thm-PPPlambdaequality}. In particular, 
Throughout the proof $X^\lambda=\{ X_i \}_{i \in I}$ denotes a homogeneous PPP with intensity $\lambda$ in $\R^d$, and we write $X^\la$ instead of $\X^\la$ since marks are non-random. Let us introduce the notion and elementary properties of Boolean models with (constant) radius $r>0$. The \emph{Poisson--Boolean model} $B(X^\lambda,r)$ (with constant connection radii $r$) is defined as
\[ B(X^\lambda,r)= \bigcup_{i \in I} B_{r}(X_i)=X^\lambda \oplus B_r(o). \]
%where for $x \in \R^d$ and $r>0$, $B_r(x)$ denotes the open Euclidean ball of radius $r$ centered at $x$. 
Connecting any two different points $X_i,X_j \in X^\lambda$ by an edge whenever 
\[ |X_i-X_j|<2r, \numberthis\label{randomGilbert} \] 
we obtain the Poisson--Gilbert graph $g_{2r}(X^\lambda)$ with connection radius $2r$. Percolation in this Gilbert graph is equivalent to the existence of an unbounded connected component in $B(X^\lambda,r)$, which we also refer to as percolation. This way, one can speak about subcritical, critical and supercritical Poisson--Boolean models.

%Now, we have the following scaling relation of the Poisson--Boolean model \cite[Section 2.2]{MR96}. For $\lambda,r>0$, $g_1(X^\lambda)$ has the same distribution as $r^{-1} g_r(X^{\lambda r^{-d}})$ and similarly, $X^\lambda \oplus B_1(o)$ has the same law as $r^{-1} (X^{\lambda r^{-d}} \oplus B_r(o))$. Consequently, percolation of $g_r(X^\lambda)$ (equivalently, $X^\lambda \oplus B_{r/2}(o)$) is equivalent to percolation of $g_{r'}(X^{\lambda'})$ (equivalently, $X^{\lambda'} \oplus B_{r'/2}(o)$) in case the following scaling relation holds:
%\[ \lambda r^d = \lambda' r'^d. \]
Recall the definition of the radius $\rG$ from~\eqref{GilbThresh} and let us fix $\la>\lc(\rG)$ for the remainder of this section. Thanks to scale invariance of Poisson--Boolean models \cite[Section 2.2]{MR96} and the well-behavedness of $\ell$, we can fix $r\in (d_o,\rG)$ such that the Poisson--Boolean model $B(X^\lambda,r/2)$ associated to $g_{r}(X^\lambda)$ is still supercritical. 
%In order to align with the notation of \cite{HJC17}, it will be useful to represent $\Xi$ as the support of the intensity measure $\Lambda$ of the modulated Poisson process given as $\Lambda (\d x) = \lambda_1 \mathds 1 \{ x \in \Xi \} \d x$, where $\lambda_1>0$ is chosen such that $\E[\Lambda(Q_1)]=1$. The modulated Poisson point process itself is the Cox point process with intensity measure $\lambda\Lambda$ for $\lambda>0$; see \cite[Example 2.1]{HJC17} for a more general formulation.  
The next lemma is an immediate consequence of the results in 
%\cite[Section 2.1]{HJC17}, which in turn follows from 
\cite[Section 1]{PP96}.
\begin{lem}[\cite{PP96}]\label{lemma-asessconn2}
Let $B(X^\lambda,r/2)$ be a supercritical Poisson--Boolean model and let $x \in \R^d$. With probability tending to one as $n \uparrow \infty$, we have that 
\begin{enumerate}[(1)]
\item\label{first-asessconn2} $B(X^\lambda,r/2) \cap Q_n(x)$ contains a connected component of diameter at least $n/3$,
\item\label{second-asessconn2} any two connected components of $B(X^\lambda,r/2) \cap Q_n(x)$ of diameter at least $n/9$ each are contained in the same connected component of $B(X^\lambda,r/2) \cap Q_{2n}(x)$.
\end{enumerate}
\end{lem}
Using Lemma~\ref{lemma-asessconn2}, we construct a renormalized percolation process on $\Z^d$. For $z \in \Z^d$, let $\Xi_n(z)$ denote the union of all connected components of $B(X^\lambda,r/2) \cap Q_n(z)$ that are of diameter at least $n/3$. For $n  \geq 1$, we say that the site $z \in \Z^d$ is \emph{$n$-good} if 
\begin{enumerate}[(1)]
\item $\Xi_n(nz) \neq \emptyset$, and
\item for any $z' \in \Z^d$ with $|z-z'|_\infty \leq 1$, it holds that all pairs of connected components $C$ of $\Xi_n(nz)$ and $C'$ of $\Xi_n(nz')$ are contained in the same connected component of $B(X^\lambda,r/2) \cap Q_{6n}(nz)$. 
\end{enumerate}
The site $z \in \Z^d$ is \emph{$n$-bad} if $z$ is not $n$-good. We have the following lemma.
\begin{lem}\label{lemma-fixedlambdagoodness}
Under the assumptions of Theorem~\ref{thm-PPPlambdaequality}, for all $n \geq 1$ sufficiently large, there exists $q_{A}=q_A(\lambda,n) \in (0,1)$ such that for any $N \in \N$ and pairwise distinct $z_1,\ldots,z_N \in \Z^d$ we have
\[ \P(z_1,\ldots,z_N\text{ are all $n$-bad}) \leq q_A^N. \]
Further, for any $\eps>0$, for all large enough $n$ one can choose $q_A$ such that $q_A<\eps$.
\end{lem}
\begin{proof}
For $z \in \Z^d$, $\mathds 1 \{ z\text{ is $n$-good} \}$ is measurable with respect to $X^{\lambda} \cap(Q_{6n}(nz) \oplus B_{r/2}(o))$, which is contained in $X^\lambda \cap Q_{7n}(nz)$ for all $n$ large enough, hence for all sufficiently large $n$ the process of $n$-good sites is 7-dependent thanks to the independence property of the PPP $X^\lambda$. Hence, using a standard argument (using dependent percolation theory~\cite{LSS97}, like in the proof of Lemma~\ref{prop-interferencecontrolIN}), it suffices to verify that 
\[ \limsup_{n\uparrow\infty}\P(o\text{ is $n$-bad})=0. \numberthis\label{badness} \]  
The limit~\eqref{badness} can be verified along the lines of the proof of~\cite[Theorem 2.6]{HJC17} using an adequate interpretation of the Poisson--Boolean model. More precisely, in view of Definition~\ref{defn-asessconn}, the assertion of Lemma~\ref{lemma-asessconn2} is equivalent to the statement \cite[Section 2.1]{HJC17} that the (for all sufficiently large $b>0$) $b$-dependent directing random measure $\Lambda$ given as $\Lambda(\d x)= \lambda_1 \mathds 1 \{ x \in B(X^\lambda,r/2) \} \d x$
is asymptotically essentially connected, where $\lambda_1>0$ is such that $\E[\Lambda(Q_1)]=1$. 
%The associated CPP is a special case of the so-called \emph{modulated PPP}, see \cite[Section 5.2.2]{CSK+13}.
%The interesting feature is that this asymptotic essential connectedness holds for any supercritical Poisson--Boolean model.
\end{proof}
The other essential proof ingredient is the interference control. We recall the ``shifted'' path-loss functions $\ell_a$ \eqref{elladef} and the shot-noise processes $I_a(x),I(x)$ from Section~\ref{sec-phasetransitionproof}, and also that by the triangle inequality, for $a \geq 0$, $I(x) \leq I_{a}(z)$ holds for any $z \in \R^d$ and $x \in Q_a(z)$.

For $n \geq 1$ and $M>0$, we say that $z \in \Z^d$ is \emph{$(n,M)$-tame} if $I_{7n}(nz) \leq M$ and \emph{$(n,M)$-wild} otherwise. Then we have the following assertion, which holds for all $\lambda$ such that $B(X^\lambda,r/2)$ is supercritical.
\begin{lem}{\cite{T18}}\label{lemma-interferencecontrol}
Under the assumptions of Theorem~\ref{thm-PPPlambdaequality}, for fixed $n \geq 1$, for all sufficiently large $M>0$, there exists $q_B=q_B(\lambda,n,M)\in (0,1)$ such that for any $N \in \N$ and pairwise distinct $z_1,\ldots,z_N \in \Z^d$ we have
\[ \P(z_1,\ldots,z_N\text{ are all $(n,M)$-wild}) \leq q_B^N. \]
Further, for $\eps>0$, for any $n \geq 1$, for all sufficiently large $M$ one can choose $q_B$ such that $q_B < \eps$.
\end{lem}
\begin{proof}
Clearly, the Lebesgue measure $\Lambda$ is asymptotically essentially connected, $b$-dependent for any $b>0$, and $\Lambda(Q_1)$ has all exponential moments. Hence the lemma can be proven very similarly to \cite[Proposition 3.3]{T18} under the condition (2b) in \cite[Theorem 2.4]{T18}. The only difference is that in \cite{T18}, $I_{6n}(nz)$ was considered instead of $I_{7n}(nz)$, but this makes no qualitative difference for the proof. Further, the additional condition in \cite[Theorem 2.4]{T18} that $\ell(0) \leq 1$ can be assumed to hold without loss of generality, for the same reason as in the proof of Proposition~\ref{prop-interferencecontrol} (see at the beginning of Step~\ref{step-interferencecontrol} in the proof of Theorem~\ref{thm-randomfadings}). 
\end{proof}
Equipped with these results, we can now prove our main theorem. 
\begin{proof}[Proof of Theorem~\ref{thm-PPPlambdaequality}.] 
For $n \geq 1$ and $M>0$, we say that the site $z \in \Z^d$ is $(n,M)$\emph{-nice} if it is both $n$-good and $(n,M)$-tame. We claim that for all sufficiently large $n$ and accordingly chosen large enough $M$, the process of $(n,M)$-nice sites percolates. Indeed, this follows by combining the estimates of Lemmas~\ref{lemma-asessconn2} and \ref{lemma-interferencecontrol} similarly to Corollary~\ref{prop-twoPeierls} and carrying out a Peierls argument. 

We claim that this assertion implies percolation in $G_{\g}(X^\lambda)$ for small $\gamma>0$. Indeed, let $n,M$ be so large that the process of $(n,M)$-nice sites percolates, and such that $Q_{6n}(o) \oplus B_{r/2}(o) \subseteq Q_{7n}(o)$. Using a standard argument (cf.~\cite{DF06} or Step~\ref{step-alsoSINRperc} in the proof of Theorem~\ref{thm-randomfadings} Part~\eqref{second-fading}), one can choose $\gamma>0$ sufficiently small such that for any $(n,M)$-tame site $z$, all connections in $g_r(X^\lambda) \cap Q_{7n}(nz)$ also exist in $G_{\g}(X^\lambda) \cap Q_{7n}(nz)$. 

Now, analogously to \cite[Section 5.2]{HJC17}, we can argue as follows. Let $\mathcal C$ be an infinite connected component of the process of sites that are $(n,M)$-nice. Let $z,z' \in \mathcal C $ and $\{ z_0=z,z_1,\ldots,z_{k-1},z_k=z' \}$ a path in $\mathcal C$ connecting $z$ and $z'$. Then, thanks to $n$-goodness, for any $j=0,\ldots,k$ and for any $X_j \in X^\lambda$ such that $B_{r/2}(X_j) \cap Q_n(nz_j) \subseteq \Xi_n(nz_j)$ we have that $X_j$ and $X_{j+1}$ are in the same connected component of $B(X^\lambda,r/2)\cap Q_{6n}(nz_j)$.
%belong to the same connected component of $B(X^\lambda,r/2) \cap Q_{6n}(nz_j)$. In other words, $X_i$ and $X_j$ are connected by a path in $g_r(X^\lambda) \cap (B_{r/2}(X_i) \cup B_{r/2}(X_j))$. 
In other words $X_j$ and $X_{j+1}$ are connected in the Poisson--Gilbert graph $g_r(X^\lambda)$ via a path in  $Q_{7n}(nz_j)$, where the additional unit of $n$ comes from the fact that centers of balls in the Boolean model might lie in a neighboring box. 
%Further, any $X \in X^\lambda$ such that $B_{r/2}(X) \cap Q_{6n}(nz'') \neq \emptyset$ for some $z''\in\Z^d$, lies in $Q_{7n}(nz'')$. 
Hence, using $(n,M)$-tameness, we conclude that all edges of this path in $g_r(X^\lambda)$ also exist in $G_{\g}(X^\lambda)$. Thus, $G_{\g}(X^\lambda)$ also percolates.
Since $\lambda>\lambda_{\rm c}(\rG)$ was arbitrary, the theorem follows.
\end{proof}

\end{document}